\documentstyle[12pt]{amsart}

\makeatletter                                            %This piece serves
\renewcommand{\@begintheorem}[2]{                        %to put numbers
\rm \trivlist \item [\hskip \labelsep {\bf #2\ \ #1.}]   %of theoremlike
                                }                        %environments
\makeatother                                             %in front.

\newcommand{\newsubsection}%
{{\bf\refstepcounter{subsection}\thesubsection\ \ }}% Makes new
\newcommand{\newsubsubsection}%
{{\bf\refstepcounter{subsubsection}\thesubsubsection\ \ }}% Makes
\newcommand{\ts}{\vspace{\baselineskip}{\bf Proof.}$\;\;$}
\newcommand{\ZZ}{{\bf Z}}
\newcommand{\QQ}{{\bf Q}}
\newcommand{\RR}{{\bf R}}

\newcommand{\CC}{{\bf C}}

\newcommand{\FF}{{\bf F}}

\newcommand{\PP}{{\bf P}}

\newcommand{\cC}{{\cal C}}

\newcommand{\cI}{{\cal I}}

\newcommand{\cM}{{\cal M}}
\newcommand{\cN}{{\cal N}}
\newcommand{\cO}{{\cal O}}

\newcommand{\cT}{{\cal T}}

\newcommand{\tT}{{\tilde{T}}}

\newcommand{\hB}{{\hat{B}}}
\newcommand{\hC}{{\hat{C}}}
\newcommand{\hT}{{\hat{T}}}
\newcommand{\hV}{{\hat{V}}}

\begin{document}

\title[The Chow group of the moduli space of marked cubic surfaces]
{The Chow group of the moduli space of marked cubic surfaces}

%\thanks{}
\author{Elisabetta Colombo}
\address{Dipartimento di Matematica, Universit\`a di Milano,
  via Saldini 50, I-20133 Milano, Italia}
\email{Elisabetta.Colombo@@mat.unimi.it}
%\thanks
\author{Bert van Geemen}
\address{Dipartimento di Matematica, Universit\`a di Milano,
  via Saldini 50, I-20133 Milano, Italia}
\email{geemen@@mat.unimi.it}
%\thanks{}

\begin{abstract}
Naruki gave an explicit construction of the moduli space of
marked cubic surfaces, starting from a toric variety and proceeding
with blow ups and contractions. Using his result, we compute
the Chow groups and the Chern classes
of this moduli space.
As an application we relate a recent result of Freitag on the
Hilbert polynomial of a certain ring of modular forms
to the Riemann-Roch theorem for the moduli space.
\end{abstract}

\maketitle

\

\begin{flushright}
{\it Dedicated to the memory of our friend Fabio Bardelli}
\end{flushright}

\

Following on the work of Allcock, Carlson and Toledo \cite{ACT},
which identified the moduli space of marked cubic surfaces $\cM$
as a ball
quotient, there has been a renewed interest in moduli spaces of
cubic surfaces. In particular, Allcock and Freitag \cite{AF}
found a projective embedding of $\cM$ using new results of
Borcherds on modular forms.
This map was actually described earlier by Coble,
who identified $\cM$ with the moduli space of six points in $\PP^2$.
We will thus call this map the CAF-map.

The moduli space $\cM$ is smooth
except for $40$ singular points, the cusps. Blowing up the cusps, one
obtains a smooth projective variety $\cC$ which we refer to as Naruki's
cross ratio variety. Using basic work of Cayley on cubic surfaces and
associated projective invariants, certain cross ratios, Naruki
\cite{Naruki}
realized that $\cC$ could be obtained from a toric variety associated
to the root system $D_4$ via a process of explicit blow ups and
contractions. Moreover he showed that there is a biregular action
of $W(E_6)$ on $\cC$.
His description of $\cC$ is at the basis of this paper.

We determine the Chow groups of $\cC$ and the spaces
of $W(E_6)$-invariant cycles in section \ref{chowgr}.
In section \ref{cusps} we consider the exceptional
divisors in $\cC$ over the cusps in $\cM$, which we call cusp divisors
(these are actually very simple varieties, being the product of three
$\PP^1$'s). We study the $36$ boundary divisors,
which parametrize nodal cubic surfaces,
in section \ref{boundary}.
The image of such a divisor in $\cM$
is isomorphic to the Segre cubic in $\PP^4$.
We use the information obtained on these divisors to compute
the Chern classes of $\cC$ in section \ref{chern}.
In section \ref{tritangent} we consider the $45$ tritangent divisors,
these parametrize cubics with an Eckardt point
(a point on the cubic surface through which 3 lines pass).
All these divisors correspond, in a $W(E_6)$-equivariant way,
to points in a finite
projective geometry. For the convenience of the reader we added some
tables at the end of the paper which give such a correspondence.

The main application of our results
is an explicit form of the Riemann-Roch theorem on
the fourfold $\cC$ in section \ref{chern}.
It allows us to compute the Euler characteristic
of integral multiples of the divisor class which gives the CAF-map.
The quartic polynomial we found (see Theorem \ref{poli})
agrees with the Hilbert function of
a certain graded ring of modular forms associated to the CAF-map
which was recently determined by Freitag \cite{Freitag}.
We also obtain some results on the Picard group of the moduli space
of smooth marked cubic surfaces in section \ref{cm0}.

Finally, we should probably point out that in using Naruki's
description, we do not (need to) consider families of
cubic surfaces and the `tautological' classes associated to them.
This is somewhat unfortunate, as this might give more insight in the
structure of the Chow groups of $\cC$ viewed as compactification of
a moduli space.

\section{The cross ratio variety}

\subsection{The toric variety \mbox{$\tT$}.}
The starting point of Naruki's construction is the smooth toric
variety $T\cong(\CC^*)^4\hookrightarrow \tT$
defined by the fan in the weight lattice
$N\cong\ZZ^4$ of the root system of type
$D_4$ whose 4-dimensional cones are the 192 Weyl chambers,
see \cite{Naruki} and \cite{vG}, section 2.
This fan has 48 edges (i.e.\ 1-dimensional cones),
which correspond to the divisors in the boundary of $\tT$.
Identifying $N\otimes_\ZZ\RR$ with $\RR^4$, with standard
basis $\epsilon_1,\ldots,\epsilon_4$ and standard scalar product, let
$$
S:=\{\pm \epsilon_i\}\,\cup\,\{(\pm
\epsilon_1\pm\epsilon_2\pm\epsilon_3\pm\epsilon_4)/2\},\qquad
R:=\{\pm\epsilon_i\pm\epsilon_j\},
$$ 
then $S$ and $R$ each have 24 elements and the half-lines
they generate are the edges of the fan.
The elements from $S\cup R$ generate the $\ZZ$-module $N$,
for $\tau\in S\cup R$ we denote by $V(\tau)\;(\subset\tT)$
the corresponding divisor.

The character group of $T$ is identified with $M:=Hom(N,\ZZ)$,
the root lattice of $D_4$. A $\ZZ$-basis of $M$ is
$$
e_1-e_2,\quad e_2-e_3,\quad e_3-e_4,\quad e_3+e_4,
$$
where the $e_i$ are the dual basis of the $\epsilon_i$,
the corresponding characters $T\rightarrow \CC^*$ are usually denoted
by $\lambda$,
$\rho$, $\nu$ and $\mu$ respectively.

\subsection{The cross ratio variety \mbox{$\cC$}.}\label{crv}
Naruki's (smooth, projective) cross ratio variety $\cC$ is obtained
from the
toric variety $\tilde{T}$ as follows (\cite{Naruki}, \S 10-12):
$$
\begin{array}{ccrcccccc}
\cC&\stackrel{r}{\longleftarrow}&
\hat{T}&\stackrel{\pi''}{\longrightarrow}&
             \tT''&\stackrel{\pi'}{\longrightarrow}&
             \tT'&\stackrel{\pi_e}{\longrightarrow}&\tT.\\
%&&r\downarrow&&&&&&\\
%\cM&\longleftarrow&\cC&&&&&
\end{array}
$$

The map $\pi_e$ is the blow up of $\tT$ in the identity element
$e\in T$.
The exceptional divisor $\pi_e^{-1}(e)\cong \PP^3$ is denoted by
$\PP^3_{\rm w}$.
The image in $\cC$ of its strict transform $\hat{E}$ in $\hT$ is a
tritangent divisor denoted by $D_{\rm w}=D_{(16)}$.

The map $\pi'$ is the blow up of $\tT'$ in $12$ smooth,
disjoint, rational curves
$C_1',\ldots,C_{12}'$.
We write $C_j'':=(\pi')^{-1}(C_j')$ for the exceptional divisor,
the map $C_j''\rightarrow C_j'$ induced by $\pi'$ is a $\PP^2$-bundle.

The map $\pi''$ is the blow up in % the strict transform
$16$ smooth disjoint surfaces
$S_1'',\ldots,S_{16}''\subset \tT''$.
Each exceptional divisor $C_j''$ meets $4$
of the surfaces $S_i''$
in $4$ disjoint rational curves.
These $4$ curves
are sections of
the $\PP^2$-bundle $C_j''\rightarrow C'_j$
 which meet
each fiber of the $\PP^2$-bundle $C_j''\rightarrow C_j'$
in $4$ general points (cf.\ \cite{Naruki}, Prop.\ 9.1, \S 10).
The strict transforms $\hat{C}_{j}$ in $\hat{T}$ of
the $C_j''$ are trivial bundles:
$$
\hat{C}_{j}\cong \PP^1\times V\,
\stackrel{\pi'\circ \pi''}{\longrightarrow}\,
\PP^1\cong C_j'
$$
here $V$ is the blow up
of $\PP^2$ in $4$ points, the map
induced by $\pi'\circ \pi''$ is the projection on the first factor.
The map $r:\hT\rightarrow \cC$ induces the projection
$\hat{C}_j\rightarrow V$ on the second factor
and $r$ is an isomorphism on the complement of the $12$ $\hat{C}_i$'s
(\cite{Naruki}, p.\ 22, 23 and Prop.\ 11.3).

The 16 exceptional divisors $\hat{S_i}$ of $\pi''$
in $\hat{T}$ map under $r$
to divisors in $\cC$.
The $W(E_6)$-orbit of these $16$ divisors in $\cC$
consists of $40$ divisors, the other 24 are
the images under $r$ of the strict transforms of the
$V(\tau)$'s with $\tau\in R$ (\cite{Naruki}, Prop.\ 11.2).
We call these 40 divisors the cusp divisors of $\cC$.

\subsection{The marked moduli space \mbox{$\cM$}}\label{cM}
There is a morphism
$$
c:\cC\longrightarrow \cM,
$$
where $\cM$ is the moduli space of
semistable marked cubic surfaces, which contracts the $40$
cusp divisors to (singular) points
(cf.\ \cite{Naruki}, Introduction and \S 12),
the cusps of $\cM$. The map $c$ is an isomorphism between
the complement of the cusp divisors in $\cC$ and the complement
of the 40 cusps of $\cM$.
The Weyl group $W(E_6)$ acts biregularly on both $\cC$ and $\cM$,
the morphism $c$ is
$W(E_6)$-equivariant.

There is a $W(E_6)$-equivariant embedding, the CAF-map
$$
F:\cM\longrightarrow \PP^9
$$
which was studied extensively by Coble, Allcock and Freitag
\cite{AF}, \cite{Freitag2} and also in \cite{vG}.

\subsection{The cusp divisors.}
The cross ratio variety $\cC$ has 40 cusp divisors of $\cC$,
each of these is isomorphic to
$\PP^1\times \PP^1\times \PP^1$,
cf.\ \cite{Naruki}, \S 12.

\subsection{The boundary divisors.}
A boundary divisor is an irreducible component of
the closure of the set of marked cubics with a node.
There are 36 such divisors in $\cC$, 24 of which are obtained from the
$V(\tau)$ with $\tau\in S$, the other 12 come from the subtori
defined by $\alpha=1$ in $T$ where $\alpha$ is a positive root of
$D_4$.

\subsection{The tritangent divisors.}
A tritangent divisor is an irreducible component of
the closure of the set of marked cubics with an Eckardt point,
that is a point on a cubic surface where three lines meet.
There are 45 such divisors in $\cC$. One of them is obtained from
$\PP^3_{\rm w}$, the other 44 are defined by the explicit equations
in \cite{Naruki}, Table 3.

\subsection{Incidence between divisors.}\label{ni}
The cusp, boundary and tritangent divisors can be very conveniently
parametrized by points in a finite projective orthogonal geometry.
Let
$$
b:\FF_3^5\times\FF_3^5\longrightarrow \FF_3,\qquad
b(x,y)=\sum_i x_iy_i
$$
be the standard bilinear form on $\FF_3^5$ and let $q(x)=b(x,x)$
be the associated quadratic form. The orthogonal group
$O(\FF_3^5,q)$ is isomorphic to $W(E_6)\times\{\pm 1\}$.
Note that in \cite{AF} the quadratic form $x_1^2-x_2^2-\ldots-x_5^2$
is used, substituting $x_2:=x_2+x_3$, $x_3=x_2-x_3$ and similarly for
$x_4,\,x_5$, one obtains $q(x)$ since $-2\equiv 1$ mod 3.

As the nonzero elements in
$\FF_3$ are $\pm 1$, the subsets
$$
N_i:=\{z\in \PP(\FF_3^5):\;q(z)=i\;\}
$$
for $i\in \FF_3$ are well defined. It is easy to check that
$$
\sharp N_0=40,\qquad \sharp N_1=45,\qquad \sharp N_{-1}=36.
$$
The group $W(E_6)$ acts transitively on each of these sets.

There is a natural $W(E_6)$-equivariant
identification between $N_0$ and the
set of cusps, $N_1$ and the set of tritangent divisors,
$N_{-1}$ and the set of boundary divisors, see \cite{AF}, section 3.
There it is also shown that a cusp divisor and a boundary
(tritangent) divisor have non-empty intersection
iff the corresponding lines in $\FF_3^5$ are perpendicular.
In particular, a boundary divisor meets 10  and
a tritangent divisor meets 16 cusp divisors.

Lemma 3.1 of \cite{AF} shows that
if $D_z$ and $D_w$ are boundary divisors corresponding to
$z,\,w\in N_{-1}$ then $D_z\cap D_{w}$ is non-empty
iff the lines $z,w\subset\FF_3^5$ are perpendicular.
This implies that a boundary divisor meets 15 other boundary divisors.

A tritangent divisor has non-empty intersection
with any other tritangent
divisor (\cite{vG}, 6.7) and with any boundary divisor (for this it
is most convenient to use the description of these divisors as images
of long and short mirrors in the 4-ball from \cite{AF}, section 3).

If $v\in \FF_3^5$ is a vector with $q(v)\neq 0$, the map
$s_v:x\mapsto x+q(v)b(x,v)v$ is in $O(\FF^5_3,q)$,
note that $s_v=s_{-v}$.
For $v\in N_{-1}$ the map $s_v$ corresponds
to the reflection in a roots of $E_6$.
For $v\in N_1$, the map $-s_v$ corresponds to an involution in
$W(E_6)$. The divisor $D_v$, for $v\in N_{-1}\cup N_1$,
is the unique divisor in the fixed point set of
the involution $\pm s_v$ in $\cC$.

\section{The Chow groups}\label{chowgr}

\subsection{} For a smooth variety $X$ we denote by $A^i(X)$ the
Chow group of codimension $i$-cycles modulo rational equivalence.
For a ring $R$ we write $A^i(X)_R:=A^i(X)\otimes_\ZZ R$.
We determine the Chow groups of $\cC$ in Theorem \ref{chowC}, some
results on $A^2(\cC)$ are actually proved in \ref{ch2i}.
Various intersection numbers are computed in \ref{canbunC} and
\ref{intcusp}.

\subsection{}
There are three divisors on $\cC$ which are obviously invariant under
$W(E_6)$, they are:
$$
\hB,\qquad\hC,\qquad\hT
$$
which are the sum of the $36$ boundary divisors,
the sum of the $40$ cusp divisors and the sum of the
$45$ tritangent divisors (actually we already used the notation
$\hT$ in \ref{crv}, however no confusion is possible).
Theorem \ref{chowC}.\ref{relT} shows
that these classes are
linearly dependent in $A^1(\cC)$.

\subsection{Symmetrizing.}\label{symm}
Let $D$ be a divisor on $\cC$ which is a linear combination
of tritangent, boundary and cusp divisors. If the class of $D$ in
the Chow group
$A^1(\cC)$ is invariant, then $\sigma^*D$ has the same class
for any $\sigma\in W(E_6)$. Therefore the class of $D$ is the
same as the class of the sum over $\sigma\in W(E_6)$ of the $\sigma^*D$
divided by $\sharp W(E_6)$.

Since $W(E_6)$ permutes the 45 tritangent divisors transitively
(and similarly the 36 boundary divisors as well as
the 40 cusp divisors) the net result is that in $D$ we replace
each tritangent divisor by $\mbox{$\frac{1}{45}$}\hT$ etc.

\subsection{Theorem.}\label{chowC}
\begin{enumerate}
\item{}
The Chow groups of the cross ratio variety $\cC$ are free
$\ZZ$-modules of rank:
$$
{\rm rk}\,A^0(\cC)={\rm rk}\,A^4(\cC)=1,\qquad
{\rm rk}\,A^1(\cC)={\rm rk}\,A^3(\cC)=61,\qquad
{\rm rk}\,A^2(\cC)=147.
$$
\item \label{relT}
The subgroups of $W(E_6)$-invariant classes have rank:
$$
{\rm rk}\,A^i(\cC)^{W(E_6)}=1,\qquad
{\rm rk}\,A^j(\cC)^{W(E_6)}=2\qquad
(i=0,\,4,\quad j=1,\,2,\,3).
$$
The classes $\hB^i,\,\hC^i$ are a basis of $A^i(\cC)^{W(E_6)}_\QQ$,
we have:
$$
\hT=(25\hB+27\hC)/4\;\in A^1(\cC),\qquad
\hB\hC=-3\hC^2\; \in A^2(\cC).
$$

\item
The canonical class of the cross ratio variety is:
$$
K_{\cC}=(-\hB+\hC)/4.
$$

\item
The decomposition of $A^1(\cC)$ and $A^3(\cC)$ into irreducible
$W(E_6)$-representations is:
$$
A^1(\cC)_\CC\cong A^3(\cC)_\CC =
{\bf 1}\oplus{\bf 1}\oplus{\bf 15_q}\oplus{\bf 20_p}\oplus{\bf 24_p}.
$$
\item
The cusp divisors are linearly independent in $A^1(\cC)$,
the $40$-dimensional subspace of $A^1(\cC_\CC)$ which they span
decomposes as:
$$
{\bf 1}\oplus {\bf 15_q}\oplus{\bf 24_p}.
$$
Similarly, the boundary divisors span a $36$-dimensional subspace
of $A^1(\cC)_\CC$ which decomposes as:
$$
{\bf 1}\oplus {\bf 15_q}\oplus{\bf 20_p},
$$
and the tritangent divisors span a $45$-dimensional subspace:
$$
{\bf 1}\oplus {\bf 20_p}\oplus{\bf 24_p}.
$$
\end{enumerate}

\ts
The construction of $\cC$ via the diagram in \ref{crv} shows that
the Chow groups of $\cC$ are free $\ZZ$-modules of finite rank
and that
$A^i(\cC)\cong H^{2i}(\cC,\ZZ)$ whereas $H^{2i+1}(\cC,\ZZ)=0$.
Since  all varieties in the diagram in \ref{crv},
are smooth, so $A^1\cong A^3$, we concentrate on $A^1$ and $A^2$.

The rank of $A^1(\tT)$ is given by the number of
edges in the fan defining $\tT$, which is $48$, minus the rank
of the torus, which is $4$. The Euler characteristic of $\tT$ is
equal to the number of maximal cones (\cite{Fulton}, p.59),
which is $192$.
Hence the rank of $A^2$ is $192-2-2\cdot 44$:
$$
{\rm rk}\,A^1(\tT)=44,\qquad {\rm rk}\,A^2(\tT)=102.
$$

The Chow groups of the blow up $\tT'$ of $\tT$ in the
identity element are $A^i(\tT')\cong A^i(\tT)\oplus A^i(\PP^3_{\rm w})$,
for $i=1,\,2$ hence the ranks increase by $1$.

The variety $\tT''$ is obtained from $\tT'$ by blowing up
the $12$ disjoint $C_j'$'s, each of which is a
$\PP^1$, hence we get:
$$
{\rm rk}\,A^1(\tT'')=45+12=57,\qquad
{\rm rk}\,A^2(\tT'')=103+24=127.
$$

Next one blows up $16$ surfaces $S_i''$ in $\tT''$ to obtain the
variety $\hT$. These surfaces are the strict transforms in $\tT''$
of surfaces $S_i\subset \tT$. Each $S_i$ is the closure of a subtorus
(for example $\nu=\rho=1$ defines an $S_i$, the first table
in \ref{cd} lists all $16$ surfaces). Using the description of
the fan defining $\tT$ one finds that ${\rm rk}A^1(S_i)=4$ for all $i$.
These surfaces all contain $e\in \tT$ hence ${\rm rk}A^1(S_i')=5$,
where $S_i'$ is the strict transform of $S_i$ in $\tT'$.
Each $C_j'$ is the strict transform of the
closure $C_j$ of a $1$-dimensional subtorus of $\tT$ (\cite{Naruki},
p.20), hence
$C_j$ is either contained in $S_i$ or meets it in $e$.
As $e$
gets blown up, an $S_i'$ either contains a $C_j'$ or they are disjoint.
Thus $\pi'':S_i''\rightarrow S_i'$ is an isomorphism and
${\rm rk}A^1(S_i'')=5$. Therefore:
$$
{\rm rk}\,A^1(\hT)=57+16=73,\qquad
{\rm rk}\,A^2(\hT)=127+5\cdot 16=207.
$$

Finally we consider the blow down $r:\hT\rightarrow \cC$,
it contracts the
12 $\hat{C}_j\cong \PP^1\times V$
in $\hT$ to the surfaces $V$ in $\cC$.
Since $V$ is a $\PP^2$ blown up in 4 points we get:
$$
{\rm rk}\,A^1(\cC)=73-12=61,\qquad
{\rm rk}\,A^2(\cC)=207-5\cdot 12=147.
$$

A relation between the invariant divisor classes can be obtained
from the divisor of the rational function $\lambda-1$ on $\cC$,
see section \ref{divlambda1}.
This divisor can be written as $(\lambda-1)_0-(\lambda-1)_\infty$
where $(\lambda-1)_0$ is the sum of one tritangent, one boundary and
4 cusp divisors, whereas $(\lambda-1)_\infty$ is the sum of
6 boundary and 9 cups divisors, one of which has multiplicity two.
Symmetrizing this relation in $A^1(\cC)$, we get:
$$
0=\frac{1}{45}\hT +
\frac{1}{36}(1 -6 )\hB+
\frac{1}{40}(4 -10 )\hC,
$$
from which we get the relation $4\hT=25\hB+27\hC$
stated in the theorem. The relation in $A^2$ will be proved in
Proposition \ref{intcusp} by restriction to a cusp divisor.

The canonical divisor of the toric variety $\tT$ is minus the sum
of the 48 divisors $V(\tau)$ (\cite{Fulton}, \S 4.3),
$\tau\in R\cup S$, half of which give
boundary divisors, the others give cups divisors. Blowing up
the origin, we must add $3\PP^3_{\rm w}$ to the strict transform
of $K_{\tT}$ to obtain $K_{\tT'}$
(\cite{GH}, Lemma in Chap.\ 1.4, p.187).
The exceptional divisors for $\pi'$
are contracted by $r$, so these do not contribute, however we must
add the 16 cups divisors which come from the blow up
of the surfaces in $\tT''$. Symmetrizing the result we get:
$$
K_{C}=\frac{3}{45}\hT+\frac{-24}{36}\hB+\frac{1}{40}(-24+16)\hC=
(-\hB+\hC)/4
$$
where we used the relation obtained before to eliminate $\hat{T}$.
A nice, explicit, anti-canonical divisor is given at the
end of \ref{divlambda}.

From the construction of $\cC$ it is clear that $A^1(\cC)$ is
generated by the tritangent, boundary and cusp divisors.
The divisor of $(\lambda-1)$ on $\cC$ contains exactly one tritangent
divisor with multiplicity one and this shows that the 36 boundary
divisors and the 40 cusp divisors suffice to generate $A^1(\cC)$.
Since rk$A^1(\cC)$=61, there must be 15 independent relations
between the boundary and cusp divisors. One such relation,
which involves both boundary and cusp divisors,
is given by the divisor of the rational function $\lambda$
(see \ref{divlambda}), others
can be obtained by applying $W(E_6)$ to this relation.

The group $W(E_6)$ permutes the boundary and the cusp divisors,
the decomposition of these permutation representations can be
found in the atlas \cite{atlas}, p.\ [26], $U_4(2)$,
see also \cite{Frame} whose notation we use.
The permutation representations on the 36 boundary components
and the 45 tritangents are denoted by $\chi^{(36)}$ and $\chi^{(45)}$
in \cite{Frame}, p.\ 100. The group $W(E_6)$ has two
permutation representations of degree 40 (see \cite{atlas}), the one
which corresponds to the 40 cusp divisors is not $\chi^{(40)}$ in
\cite{Frame} but it is the one listed in the theorem, as can be
verified by computing a few traces.

The only $15$-dimensional $W(E_6)$ representation which
the two permutation
representations have in common is the irreducible representation
${\bf 15_q}$, hence this must be the representation on the relations.
Therefore the representation on $A^1(\cC)$ is the sum of the
two permutation representations modulo ${\bf 15_q}$.
In particular $A^1(\cC)^{W(E_6)}$
has rank 2 and is generated (over $\QQ$) by $\hB$ and $\hC$.
Since the permutation representation on the tritangents does not
contain ${\bf 15_q}$ and since the relation provided by $\lambda-1$
involves both boundary and cusp divisors, the subrepresentations
${\bf 20_p}$ (which is in
common with the boundary divisors) and ${\bf 24_p}$ (which is in
common with the cusps) cannot map to zero in $A^1(\cC)$,
and finally since $\hT\in A^1(\cC)$ is non-zero as well,
we conclude that the tritangents divisors are independent in $A^1(\cC)$.

The proof of ${\rm rk}A^2(\cC)^{W(E_6)}=2$ will be given in \ref{ch2i}.
\qed

\subsection{Remark.} \label{cm0}
Let $\cM^0\subset \cM$ be the moduli space of smooth marked cubic
surfaces, so $\cM^0\cong \cC-{\rm Support}(\hB+\hC)$, and it is also
the complement in $\cM$ of the $36$ boundary divisors.
As we observed in the proof of the theorem, the first
Chow group of $\cC$ is generated by the classes of the boundary and
cusp divisors. In particular,
$$
A^1(\cM^0)=Pic(\cM^0)=0.
$$
(Note that relations $\hT=(25\hB+27\hC)/4$ and $K_{\cC}=(-\hB+\cC)/4$
only show that $4\hT$ and $4K_{\cC}$ restrict to zero on $\cM^0$).
Since finite cyclic subgroups of $Pic(\cM^0)$ correspond to
finite cyclic unramified coverings of $\cM^0$,
we conclude that the
maximal abelian quotient of $\pi_1(\cM^0)$ is trivial.
See \cite{Eduard} for interesting results on this fundamental
group.

\subsection{The hyperplane class.}\label{classH}
We determine $(F\circ c)^*\cO(1)$ where (cf.\ section \ref{cM})
$$
\cC\stackrel{c}{\longrightarrow} \cM\stackrel{F}{\longrightarrow}\PP^9.
$$
From the work of Allcock and Freitag \cite{AF} we know that, modulo
cusp divisors, it is given by a so called cross divisor, that is
$$
(F\circ c)^*\cO(1)\equiv D_\alpha+D_\beta+D_\gamma+D_\delta+D_t
\qquad  (\text{modulo cusp divisors}),
$$
where $\alpha,\ldots,\delta\in N_{-1}$ and $t\in N_1$ are pairwise
perpendicular. In particular, $D_\alpha,\ldots,D_\delta$ are
boundary divisors and $D_t$ is a tritangent divisor.

The divisor class $(F\circ c)^*\cO(1)\in A^1(\cC)$ is
$W(E_6)$-invariant (since $F\circ c$ is $W(E_6)$-equivariant), and
its expression in terms of invariant classes is given in the following
proposition.

\subsection{Proposition.}
The $W(E_6)$-invariant hyperplane class of the CAF map
is given by:
$$
(F\circ c)^*\cO(1)=(\hB+3\hC)/4.
$$

\ts
The contribution of the cusp divisors to
$D_\alpha+D_\beta+D_\gamma+D_\delta+D_t$ is the sum of the $24$
cusp divisors $D_n$, each with multiplicity one, where $n\in N_0$
is perpendicular to at least one of element of the set
$\{\alpha,\beta,\gamma,\delta,t\}$. This can be seen from
a local computation using the explicit expression of the coordinate
functions of $(F\circ c)$
(cf.\ \cite{vG}, 5.5 and the proof of Theorem 5.7).

Symmetrizing w.r.t.\ $W(E_6)$,  as explained in \ref{symm},
we obtain the result.
\qed

\section{Cusp divisors.}\label{cusps}

\subsection{}\label{cuspp13}
According to Naruki, \cite{Naruki}, each cusp divisor is isomorphic
to $\PP^1\times\PP^1\times\PP^1$.
In the following proposition we prove the
relation between $W(E_6)$-invariant codimension two classes
(already stated in Theorem \ref{chowC}) and we
determine some intersection numbers involving cusp divisors.

\subsection{Proposition.}\label{intcusp}\label{c03}\label{c0^2}
Let $C_0\cong (\PP^1)^3$ be a cusp divisor,
and let
$$
D_1=\{0\}\times(\PP^1)^2, \qquad
D_2=\PP^1\times\{0\}\times\PP^1,\qquad
D_3=(\PP^1)^2\times\{0\}\qquad(\in A^1(C_0)).
$$
There are nine boundary divisors which have non-empty intersection
with $C_0$. The class of the intersection of each of these
with $C_0$ is one of the $D_i$ and
$$
C_0 \hB=3(D_1+D_2+D_3).
$$

Let $B_i$ be a boundary divisor such that $B_iC_0=D_i$.
Then, for $i\neq j$:
$$
C_0B_1B_2B_3=1,\qquad C_0B_i^2B_j=0,\qquad C_0B_i^3=0,
\qquad C_0^2B_iB_j=-1,\qquad C_0^2B_i^2=0.
$$
Moreover,
$$
C_0^2=-C_0(B_1+B_2+B_3),\qquad C_0^3B_i=2,\qquad C_0^4=-6.
$$
Finally we have:
$$
\hB\hC=-3\hC^2\; \in A^2(\cC),\qquad \hC^4=-240.
$$

\ts
The toric subvariety $V(\epsilon_1+\epsilon_3)$ in $\tT$
is defined by the fan consisting of the images in
$\RR^4/\langle\epsilon_1+\epsilon_3\rangle$ of all the cones
containing $\epsilon_1+\epsilon_3$. One finds that this fan has six
edges and that $V(\epsilon_1+\epsilon_3)\cong (\PP^1)^3$.
The divisor $V(\epsilon_1+\epsilon_3)$ is birationally
isomorphic to a cusp divisor $C_0$
in $\cC$ via the maps in the diagram in \ref{crv} above.
One checks that each of these maps is
an isomorphism on this divisor, hence $C_0$ is isomorphic to
$(\PP^1)^3$.

The six edges of the fan defining $C_0=V(\epsilon_1+\epsilon_3)$
define the classes $D_i$ (each twice) in $A^1(C_0)$.
On the other hand, the edges
correspond to non-empty intersections
$V(\epsilon_1+\epsilon_3)\cap V(\tau)$.
One verifies that $\tau\in S$, hence
the $V(\tau)$ define boundary divisors. In particular, each
$V(\epsilon_1+\epsilon_3)\cap V(\tau)$ is a $D_i$.
That there are nine boundary divisors meeting $C_0$ follows from
\ref{ni}. In the notation of \ref{cd}, $V(\epsilon_1+\epsilon_3)$
is labelled as $[24.35.16]$. The nine roots in this set correspond to
the nine boundary divisors meeting $C_0$. These nine roots
are in three orthogonal $A_2$'s, the three divisors corresponding
to the roots from one $A_2$ intersect $C_0$ in the same class
in $A^1(C_0)$ because the corresponding divisors do not meet on $\cC$.
In particular $C_0\hB=3(D_1+D_2+D_3)$.

Now it follows that $C_0B_1B_2B_3=1$ because the $D_i$
intersect transversely in one point on $C_0$, for
$C_0^2B_iB_j$ see below. Finally
$C_0B_i^2B_j=C_0B_i^3=C_0^2B_i^2=0$ since $D_i^2=0$ in $(\PP^1)^3$.

To find the selfintersection $C_0^2$,
we restrict the divisor of the rational function $\lambda$
(cf.\ \ref{divlambda}) to $C_0=V(\epsilon_1+\epsilon_3)=D_v$
with $v=(1,-1,0,-1,0)\in N_0$.
Using \ref{ni} and the tables in the Appendix
to eliminate empty intersections, we find:
$C_0(C_0+B_1+B_2+B_3)=0$ where
the $B_i$ are the $V(\tau)$ with
$\tau=\epsilon_1$, $(\epsilon_1-\epsilon_2+\epsilon_3+\epsilon_4)/2$
,$(\epsilon_1-\epsilon_2+\epsilon_3-\epsilon_4)/2$,
which correspond to the roots
$h_{345},\,h_{136}$ and $h_{246}$. These three roots are
in three distinct $A_2$'s in $[24.35.16]$ so $C_0^2=-(D_1+D_2+D_3)$.
Using this relation twice we get:
$$
C_0^3B_i=C_0^2(-B_1-B_2-B_3)B_i=C_0(-B_1-B_2-B_3)^2B_i=2.
$$
Therefore $C_0^4=C_0^2\cdot C_0(-B_1-B_2-B_3)=-3\cdot 2=-6$.
We also get $C_0^2B_iB_j=-(D_1+D_2+D_3)D_iD_j=-1$
if $i\neq j$.

Since $\hC$ is the sum of $40$ disjoint cusp divisors, it follows
that $\hC^4=40\cdot -6=-240$. Finally we observe that
$C_0\hB=3(D_1+D_2+D_3)=-3C_0^2$, hence, summing over the $40$ cusps,
we get $\hC\hB=-3\hC^2$.
\qed

\subsection{Remark.}\label{120indep}
Recall from Theorem \ref{chowC} that $rk(A^2(\cC))=147$.
We identify $120$ independent codimension two cycles.

Let $C_0\cong (\PP^1)^3$ be a cusp divisor and let $D_i=C_0B_i\in C_0$
be the three divisors on $C_0$ as in Proposition \ref{intcusp}.
Then $D_iD_j=-1$, if $i\neq j$, and $D_i^2=0$ in $A^4(\cC)$.
Therefore the $D_i$ are independent in $A^2(\cC)$.
As the $40$ cusps are disjoint we get
$3\cdot 40=120$ independent codimension two classes on
$\cC$.

\section{The boundary divisors}\label{boundary}

\subsection{} In this section we study the boundary divisors
and we compute some intersection numbers on $\cC$ in Theorem
\ref{canbunC}.

\subsection{The Chow group of a boundary divisor.}
\label{chowB0}
Let $B_0$ be a boundary divisor on $\cC$.
Recall that the symmetric group $S_6$ acts on $B_0$.
The divisor $B_0$ meets
15 other boundary divisors, 10 cusp divisors and the
45 tritangent divisors. We denote the sum of these by
$$
\hB_b,\qquad \hC_b,\qquad \hT_b.
$$
Then in $A^1(B_0)$ we have:
$$
\hB_{|B_0}=\hB_b+B_0^2,
\qquad \hC_{|B_0}=\hC_b,\qquad \hT_{|B_0}=\hT_b.
$$

\subsection{Proposition.}\label{intB0}
The Chow groups of a boundary divisor $B_0$ are:
$$
A^1(B_0)\cong A^2(B_0)\cong \ZZ^{16}.
$$
We also have:
$$
A^1(B_0)^{S_6}\cong A^2(B_0)^{S_6}\cong \ZZ^{2},
$$
$\hB_b^i,\;\hC_b^i$ are a basis of $A^i(B_0)^{S_6}_\QQ$
for $i=1,\,2$ and $\hB_b\hC_b=-3\hC_b^2$.
The intersection numbers between invariant classes are:
$$
\hB^3_b=-165=-3\cdot 5\cdot 11,  \qquad
\hB^2_b\hC_b=180=2^2\cdot 3^2\cdot 5 ,
$$
$$
\hB_b\hC^2_b=-60=-2^2\cdot 3\cdot 5 ,\qquad
\hC^3_b=20= 2^2\cdot 5.
$$
Moreover we have:
$$
c_1(N_{B_0/\cC})=B_0^2=-(\hB_b+3\hC_b)/5,\qquad
K_{B_0}=-(2\hB_b+\hC_b)/5.
$$

\ts
To compute the Chow groups of a boundary divisor, we
choose $B_0$ to be the strict transform of the divisor
$V(\epsilon_1)$ on the toric variety $\tilde{T}$.
The divisor $V(\epsilon_1)$
is the toric variety defined by the fan consisting of
the images in $\RR^4/<\epsilon_1>$ of all the
cones which contain $\epsilon_1$.
One finds that $V(\epsilon_1)$ is smooth
and that the fan has 14 edges.
Thus ${\rm rk}A^1V(\epsilon_1)=14-3=11$. This divisor does not
contain $e\in \tT$, but it meets (transversally in one point)
exactly
one of the $12$ curves which get blown up in $\tT'$, it is
$C_1':\;\mu=\nu=\rho=1$, hence ${\rm rk}A^1(V(\epsilon_1)'')=12$.
The divisor $V(\epsilon_1)''$ in $\tT''$ meets exactly 4 of the 16
surfaces which get blown up, it meets of them in a
smooth rational curve.
Therefore ${\rm rk}A^1\widehat{V(\epsilon_1)}=12+4=16$.
Since $\widehat{V(\epsilon_1)}$ intersects just one of the $\hat{C}_i
\cong \PP^1\times S$
(it is $\hat{C}_1$) in a copy of $V$ (\cite{Naruki}, p.\ 22)
and $r$ contracts $\hat{C}_1$ to $V$, we conclude that
the blow down map $r:\hT\rightarrow\cC$ induces an isomorphism
of smooth 3-folds
$\widehat{V(\epsilon_1)}\cong B_0$. Hence $A^1(B_0)$ and $A^2(B_0)$
have rank 16.

Note that $K_{V(\varepsilon_1)}$
is minus the sum of the divisors corresponding to the 14 rays of
the fan (the strict transform of 8 of these divisors
is a boundary divisor, the other 6 are cusp divisors).
From the formula
for the canonical divisor of a blow up we find that
$K_{V(\epsilon_1)''}$
is obtained from the strict transform
of $K_{V(\varepsilon_1)}$ by adding twice
the class of the (irreducible) exceptional divisor, this exceptional
divisor is a boundary divisor.
Finally $K_{\widehat{V(\epsilon_1)}}=K_{B_0}$
is obtained by adding 4 cusp divisors to this.
Since $K_{B_0}$ is $S_6$-invariant, we can symmetrize w.r.t. $S_6$:
$$
K_{B_0}=\frac{1}{15}(-8+2)\widehat{B}_b+\frac{1}{10}(-6+4)\widehat{C}_b
=-(2\hB_b+\hC_b)/5.
$$

The selfintersection $B_0^2$ can be computed
from the restriction of the divisor of $\lambda$ on $\cC$
(cf.\ \ref{divlambda}) to $B_0$.
As $B_0=r_*\widehat{V(\epsilon_1)}$ occurs in the divisor
of zeroes with multiplicity one,
the relation $0=B_0\cdot (\lambda)$ implies that $B_0^2$
is a linear combination of $4$ boundary divisors in $B_0$
with coefficient $-1$ and $1$ boundary divisor with
coefficient $1$, as well as $4$ cusp divisors in $B_0$
with coefficient $-1$ and $1$ with coefficient $-2$.
Symmetrizing w.r.t.\ $S_6$ we get:
$$
B_0^2=-\frac{1}{15}(4-1)\widehat{B}_b-\frac{1}{10}(4+2)\widehat{C}_b
= -(\hB_b+3\hC_b)/5.
$$

It is amusing to check that these expressions for $K_\cC$, $K_{B_0}$
and $B_0^2$ do indeed satisfy the adjunction formula:
$$
K_{B_0}=(K_{\cC}+B_0)_{|B_0}.
$$
See \ref{invB0} for the $S_6$-invariant classes and
\ref{intbv} for the intersection numbers.
\qed

\

\subsection{The Segre cubic.}
The Segre cubic is the cubic threefold $S$ in $\PP^5$ defined by the
equations:
$$
S:\qquad
x_1+\ldots+x_6=x_1^3+\ldots+x_6^3=0.
$$
It is well know that any boundary divisor is isomorphic to the
blow up of $S$ in its 10 nodes
(the points with three coordinates $+1$ and
three coordinates $-1$).
In fact, the image of $B_0$ in $\cM\subset\PP^9$
is the Segre cubic (cf.\ \cite{vG}, 6.2). The 10 nodes are the images of
the intersection of $B_0$ with 10 cusp divisors, each such intersection
is isomorphic to the exceptional divisor over a node, thus to
$\PP^1\times\PP^1$.
The Segre cubic is obviously invariant under permutations of the $x_i$,
this action of $S_6$ is induced by the action of a
subgroup in $W(E_6)$ consisting of all the elements which fix $B_0$.

\subsection{Boundary surfaces in $S$.}
The intersection of $B_0$ with a boundary divisor $B_1$ is,
if non-empty,
isomorphic to $V$, a $\PP^2$ blown up in $4$ general points.
Thus $V$ maps to subvariety isomorphic to $\PP^2$ in $S$ and this
subvariety contains 4 of the singular points.

The divisor $B_1$ is contained in the fixed point set of an involution
in $W(E_6)$
(a reflection in a root of $E_6)$) and thus $B_0\cap B_1$ is
contained in the fixed
point set of an involution in $S_6$. Considering
the involutions in $S_6$, we find that only the products of
three commuting transpositions have a fixed point locus isomorphic to
$\PP^2$ and this $\PP^2$ is linearly embedded. For
$\sigma=(ij)(kl)(mn)\in S_6$ we have the fixed point set is
$$
\PP^2_\sigma=\{\,x_i+x_j=x_k+x_l=x_m+x_n=0\,\}\quad\subset S.
$$

\subsection{Parametrizing $S$.}\label{bbu}
It is well-known (cf.\ \cite{Hunt}, Theorem 3.1.2)
that the blow up $B_0$ of the Segre threefold in the 10 nodes
can be obtained as follows:
$$
S\stackrel{}{\longleftarrow}B_0
\stackrel{\rho_2}{\longrightarrow}
B_5\PP^3
\stackrel{\rho_1}{\longrightarrow}
\PP^3,
$$
here $\rho_1$ is the blow up of $\PP^3$ in 5 general points
$p_1,\ldots,p_5$ and $\rho_2 $ is
the blow up of $B_5\PP^3$ in the strict transforms of the 10 lines
spanned by pairs of the points.
Note that this implies that ${\rm rk}A^1(B_0)=1+5+10=16$,
in agreement with Prop.\ \ref{intB0}.
The rational map $\PP^3-\,\rightarrow S$ is given by the linear
system of quadrics which contain the 5 points.

\subsection{The Chow group of $B_0$ (bis).}
We denote the exceptional fibers of $\rho_1$ by $V'_1,\ldots,V_5'$
and their strict transforms in $B_0$ by $V_1,\ldots,V_5$.
The divisor $V_1'\cong \PP^2$
intersects the strict transform in $B_5\PP^3$ of the
$4$ lines $\langle p_1,p_i\rangle$ ($i=2,\ldots,5$)
in $4$ points in general position. Thus each $V_i$ is isomorphic to $V$.

Let $V''_{ij}=\langle p_k,p_l,p_m\rangle$ be the plane in $\PP^3$
spanned by the points with complementary indices, so $\{i,\ldots,m\}
=\{1,\ldots,5\}$. Its strict transform $V'_{ij}$ in $B_5\PP^3$
has divisor class $H'-V_k'-V'_l-V'_m$, here $H'$ is the pull-back
of $\cO(1)$ on $\PP^3$. Its strict transform $V_{ij}$
in $B_0$ has class
$$
V_{ij}=H-V_k-V_l-V_m-C_{kl}-C_{km}-C_{lm},
$$
where $C_{ij}$ is the exceptional
divisor in $B_0$ over the strict transform of the
line $\langle p_i,p_j\rangle$ (note that $V''_{ij}$ contains three of
the $10$ lines) and $H$ is the pull-back of $\cO(1)$ to $B_0$.
The intersection point $V_{ij}''$ with the line
$\langle p_i,p_j\rangle$
is blown up by $\rho_2$, so $V_{ij}$ is also isomorphic to $V$.

The divisors $H$, $V_1,\ldots V_5$, and the $10$ $C_{ij}$
are a $\ZZ$-basis of $A^1(B_0)$,
so it is easy to verify that the $5$ $V_i$'s and
the $10$ $V_{ij}$'s are independent. It will be convenient to define
$V_{i6}:=V_i$, so we can talk simply about the $15$ $V_{ij}$.

Using explicit formulas for the quadrics on the 5 points
(cf.\ \cite{Hunt}, proof of Thm.\ 3.2.1), one verifies that the
$15$ $V_{ij}$ map to the 15 $\PP_\sigma$. Thus we conclude,
for example, that in $A^1(B_0)$ we have:
$$
\hB_b=\sum_{ij}^6 V_{ij}=10H-5\hV-3\hC_b,\qquad\text{with}\quad
\hV:=\sum_{i=1}^5 V_i,\quad \hC_b=\sum_{kl}^5C_{kl}.
$$

\subsection{Invariant classes in $A^1(B_0)$.}\label{invB0}
The group $S_6$ permutes the $\PP_\sigma$ and thus the
$15$ $V_{ij}$'s transitively. Hence there is only a one
dimensional $S_6$-invariant
subspace in the subspace (of codimension $1$) of $A^1(B_0)_\QQ$
which they span.
Therefore $\dim A^1(B_0)_\QQ^{S_6}\leq 2$.
Since also $\hC_b$ is invariant and does not lie in the span of
the $V_{ij}$'s, we have $\dim A^1(B_0)_\QQ^{S_6}=2$ and
$\hB_b$, $\hC_b$ are a basis of the invariant classes.

\subsection{The Picard group of $V$}
Each of the $15$ boundary surfaces $V_{ij}$ of $B_0$ is isomorphic
to $V$ and it will be convenient to restrict divisor classes
on $B_0$ (and on $\cC$) to a $V_{ij}$ in order to compute intersection
numbers, Chern classes etc.

Since $V$ is isomorphic to the blow up of
$\PP^2$ in $4$ points, we have
$$
A^1(V)\cong \ZZ l \oplus\ZZ e_1\ldots\oplus\ZZ e_4
$$
with intersection product $l^2=1$, $e_i^2=-1$,
$l\cdot e_i=e_i\cdot e_j=0$ if $1\leq i,j\leq 4$ and $i\neq j$.
Moreover, with $c:=e_1+\ldots+e_4$:
$$
c_1(V)=-K_V=3l-c,\qquad c_2(V)=7,
$$
the formula for $K_V$ is standard (\cite{HAG}, V.3.3),
for $c_2$
one can use $\chi(\cO_V)=1$
and the Noether formula (\cite{HAG}, Appendix A, 4.1.2)
$\chi(\cO_V)=(K^2+c_2)/12$.

\subsection{Lemma.}\label{resV}
Let $V_{16}$ $(\subset B_0)$ be the strict transform
of the exceptional divisor $V_1'\;(\subset B_5\PP^3)$
over $p_1\in \PP^3$.
The restriction map on divisors is given  by:
$$
A^1(B_0)\longrightarrow A^1(V_{16}),\qquad
\hB_{b|V_{16}}=5l-3c,\qquad \hC_{b|V_{16}}=c,\qquad V^2_{16}=-l
\qquad\qquad(\in A^1(V_{16})).
$$
In particular, the restriction map
is injective on $S_6$-invariants.

\ts
Since $V_{16}$ maps to a point in $\PP^3$, the class
$H$ maps to zero in $A^1(V_{16})$.
Moreover, $V_{16}$ does not intersect the divisors
$V_{i6}$ ($i=2,\ldots,5$) nor the $4$ divisors
$V_{1k}$ with $k=2,\ldots,5$
nor the divisors $C_{kl}$ with
$k,\,l\neq 1$, thus also these map to zero.
The intersection $V_{16}$ and one of the $6$ divisors
$V_{kl}$, $2\leq k,l\leq 5$,
is transversal along a
divisor with class $l-e_i-e_j$ where $\{i,j,k,l\}=\{2,\ldots,5\}$.
Obviously $C_{1j}$
maps to $e_{j-1}$, $j=2,\ldots, 5$ (up to permutation of indices).

It remains to determine the image of $V_{16}\in A^1(B_0)$,
we write $V^2_{16}$ for this class.
Using $K_{B_0}=(-2/5)\hB_b+(1/5)\hC_b$, which restricts to
$(-2/5)(V_{16}^2+6l-2c))+(1/5)c$, and the adjunction formula :
$$
K_{V_{16}}={K_{B_0}}_{|V_{16}}+V^2_{16}\qquad
\text{we get}\quad
-3l+c=\frac{3}{5}V^2_{16}-\frac{12}{5}l+c,
$$
hence $V^2_{16}=-l$.
\qed

\subsection{Intersection numbers on $B_0$.}\label{intbv}
It is now easy to compute the intersection numbers.
If $D_1$ and $D_2$ are $S_6$-invariant divisors on $B_0$, then
$D_1D_2\hB_b=15D_1D_2V_{16}$ since $\hB_b$ is the sum of 15 surfaces
in the $S_6$-orbit of $V:=V_{16}$.
Recall that we may compute the intersection product $D_1D_2V$ on $B_0$
simply as ${D_1}_{|V}\cdot {D_2}_{|V}$ on $V$, in fact:
$$
D_1\cdot D_2\cdot V=D_1\cdot D_2\cdot j_*1_V=
j_*(j^*(D_1\cdot D_2)\cdot 1_V)=j_*(j^*D_1\cdot j^*D_2)
\quad \in A^3(B_0)
$$
were
$j:V\hookrightarrow B_0$ is the inclusion of $V$ in $B_0$ and
$1_V\in A^0(V)$ is the class of $V$, $1_V=[V]$, and
we used the projection formula and the fact that $j^*$ is a ring
homomorphism. Thus we find:
$$
\hB_b^3=15\hB_{b|V}\hB_{b|V}=15(5l-3c)^2=15(25+9\cdot(- 4))=165,
$$
Similarly, we get
$\hB_b^2\hC_b=15(5l-3c)c=180$ and $\hB_b\hC_b^2=15c^2=60$.
It remains to compute $\hC_b^3$, as $\hC_b$ is the sum of 10 disjoint
divisors in one $S_6$-orbit, this is $10C_b^3$ where $C_b$ is a cusp
divisor on $B_0$. But $C_b^3=C_0^3B_0=2$ by Proposition \ref{intcusp}.

\subsection{Intersection numbers on $\cC$.}
Using the results on the restriction to subvarieties
we can now determine some intersection numbers on the cross ratio
variety as well as its canonical class.

\subsection{Theorem.}\label{canbunC}
The intersection numbers between invariant classes on the cross ratio
variety $\cC$ are:
$$
\hB^4=-12528=-2^4 3^3 29,  \qquad
\hB^3\hC= 6480=2^4 3^4 5 ,\qquad
\hB^2\hC^2=-2^4 3^3 5 ,
$$
$$
\hB\hC^3=720=2^4 3^2 5 ,\qquad
\hC^4=-240=-2^4 3\cdot 5.
$$
Moreover, for boundary divisors $B_0,\ldots, B_3$ which have a
non-empty intersection (that is, the corresponding lines in
$\FF_3^5$ are perpendicular) we have:
$$
B_0^4=-3,\qquad B_0^3B_1=1,\qquad B_0^2B_1^2=1,\qquad
B_0^2B_1B_2=-1,\qquad B_0B_1B_2B_3=1.
$$

\ts
We follow the strategy of section \ref{intbv},
replacing
the role of the surface $V$ with that of the boundary divisor $B_0$.
Since $\hB$ is the sum of the 36 divisors in the $W(E_6)$-orbit of
$B_0$ we get, for $W(E_6)$-invariant divisors $D_1,D_2,D_3$:
$$
D_1D_2D_3\hB=36D_1D_2D_3B_0=36(D_1D_2D_3)_{|B_0}=
36({D_1}_{|B_0})({D_2}_{|B_0})({D_3}_{|B_0}).
$$
Using that $\hB_{|B_0}=\hB_b+B_0^2=(4\hB_b-3\hC_b)/5$
(cf.\ \ref{chowB0} and Prop.\ \ref{intB0}) we obtain:
$$
\hB^4=36{\hB}_{|B_0}^3=\mbox{$\frac{36}{125}$}(4\hB_b-3\hC_b)^3
=\mbox{$\frac{36}{125}$}(64\cdot(-165)-\ldots)=-2^4 3^3 29,
$$
similarly one finds:
$$
\hB^3\hC=\mbox{$\frac{36}{25}$}(4\hB_b-3\hC_b)^2\hC_b=
6480,\qquad
\hB^2\hC^2=\mbox{$\frac{36}{5}$}(4\hB_b-3\hC_b)\hC_b^2=-2160,
$$
and $\hB\hC^3=36\hC_b^3=720$.
For $\hC^4$ see Proposition \ref{intcusp}.

The intersection numbers of the boundary divisors
are easily determined by restriction to
$B_0$ or $B_0B_1=V_{ij}$ (for some $i,j$),
the restriction map $A^1(B_0)\rightarrow A^1(V_{ij})$
is given in \ref{resV}. Hence:
$$
B_0^4=N_{B_0/\cC}^3=-\mbox{$\frac{1}{125}$}(\hB_b+3\hC_b)^3=-3,
\qquad
B_0^3B_1=N_{B_0/\cC}^2\cdot V_{ij}=\mbox{$\frac{1}{25}$}(5l-3c+3c)^2=1.
$$
Similarly one computes the other intersection numbers.
\qed

\subsection{} It is now easy to compute the degree of the moduli space
of marked cubics embedded in $\PP^9$ by the CAF-system.
This computation was already done in \cite{Freitag}.
Freitag used two skew $\PP^4$'s in $\PP^9$
(eigenspaces of a subgroup of index two of $W(E_6)$)
and showed, using explicit generators of the
ideal of $\cM\subset \PP^9$ and computer algebra,
 that the projection from one to the
other induces a morphism $\cM\rightarrow \PP^4$ of degree $27$.

\subsection{Corollary.} The degree of the image in $\PP^9$ of the
moduli space $\cM$ under the CAF-map is $27$.

\ts
From Theorem \ref{classH} it follows that the degree is given by
$$
\left((\hB+3\hC)/4\right)^4=\mbox{$\frac{1}{256}$}
(\hB^4+12\hB^3\hC+\ldots+81\hC^4)=27,
$$
where we used the intersection numbers from Theorem \ref{canbunC}.
\qed

\subsection{The $W(E_6)$-ivariants in $A^2(\cC)$.}\label{ch2i}
We already observed in \ref{120indep} that there are 120 independent
codimension 2 classes (obtained as $B_iC_j$)
supported on the cusp divisors. It is easy to see
that $W(E_6)$ permutes these classes transitively, and hence
the subspace of invariant classes has dimension 1.

Next we
consider the orthogonal complement (w.r.t.\ the intersection form)
in $A^2$, it has dimension $147-120=27$. Classes in the complement
are given for example by the selfintersections $B_i^2$
of the $36$ boundary divisors (indeed, $B_iC_jB_k^2=0$, for any indices,
by Proposition \ref{intcusp}). A computation shows that the rank of the
intersection matrix $(\,B_i^2B_j^2\,)$ is $21$. The group $W(E_6)$
permutes these classes and since the permutation character is
${\bf 1}\oplus{\bf 15_q}\oplus{\bf 20_p}$, we conclude that there
is in $A^2(\cC)_\QQ$
a $21$-dimensional subspace, perpendicular to the $120$ subspace,
on which $W(E_6)$ acts via ${\bf 1}\oplus{\bf 20_p}$. In particular
we found another invariant.

To get our hands on the remaining $6$-dimensional subspace,
we considered the image of the $306$-dimensional
$\QQ$-vectorspace $W$ generated by the
$270=36\cdot 15/2$ (non-empty) intersections $B_iB_j$
of boundary divisors and the $36$ $B_i^2$.
The group $W(E_6)$ acts on $W$ by permuting the basis vectors.
 Using the previous
theorem (and a computer) we found that
the matrix of intersection products
has rank $147$, hence these classes span $A^2(\cC)_\QQ$ (the $270$
classes $B_iB_j$ span a subspace of codimension $1$).
Next one computes the kernel of $W\rightarrow A^2(\cC)_\QQ$
(which is the kernel of the intersection matrix). Its perpendicular
w.r.t.\ to the standard innerproduct on $W$ is a $W(E_6)$-invariant
subspace which maps isomorphically onto $A^2(\cC)_\QQ$. Now one
computes the $6$-dimensional subspace which is perpendicular to the
$120+21=141$-dimensional subspace we already found and one verifies
that it is the standard $6$-dimensional representation ${\bf 6_p}$
of $W(E_6)$.

A class in this subspace, which corresponds to the root
$h$, is the class $S_h=S_h^+ - S_h^-$, where $S_h^+$ is the sum of
the $15\cdot 4=60$ classes $B_{ij}\cdot B_{ijk}$ ($1\leq i\le j\leq 6$
and $k\neq i,j$) and $S_h^-$ is the sum of
the $60$ classes $B_{ij}\cdot B_{klm}$ ($1\leq i\le j\leq 6$
and $\{k,l,m\}\cap \{ i,j\}=\emptyset$), here we wrote $B_m$ for
the boundary divisor corresponding to the root $h_m$ (cf.\ \ref{bd}).

In any case, we can now conclude that ${\rm rk} A^2(\cC)^{W(E_6)}=2$.
As a byproduct, we found that the complement of the $120$ classes
supported on the cusps has the representation
${\bf 1}\oplus{\bf 20_p}\oplus{\bf 6_p}$, which is the permutation
representation of $W(E_6)$ on the 27 lines (!),
cf.\ \cite{Frame}, (5.2).
For completeness sake, the representation on the $120$-dimensional
subspace is
${\bf 1}\oplus{\bf 15_q}\oplus{\bf 20_p}
\oplus{\bf 24_p}\oplus{\bf 60_p}$.

\section{The Riemann-Roch Theorem}\label{chern}
\subsection{Riemann-Roch.}
The Riemann-Roch theorem for a line bundle $\cO_\cC(D)$ on the
cross ratio variety $\cC$ is:
$$
\chi(\cO_\cC(D))={\rm deg}\left(
(1+D+\mbox{$\frac{1}{2}$}D^2+\mbox{$\frac{1}{6}$}D^3+
\mbox{$\frac{1}{24}$}D^4)\cdot
(1+\mbox{$\frac{1}{2}$}c_1+\mbox{$\frac{1}{12}$}(c_1^2+c_2)+
\mbox{$\frac{1}{24}$}c_1c_2+td(\cT)_4)\right)_4
$$
where $td(\cT)_4
=\mbox{$\frac{-1}{720}$}(c_1^4-4c_1^2c_2-3c_2^2-c_1c_3+c_4)$,
(\cite{HAG}, Appendix A.4) and $\cT$ is the tangent
bundle of $\cC$. In particular, if $D=0$ we find
$1=\chi(\cO_\cC)=td(\cT)_4$. Hence to use the Riemann-Roch
theorem we only need to compute
$c_1$ and $c_2$ of $\cC$, and we will only indicate how we computed
$c_3$ and $c_4$.

\subsection{Chern classes.}\label{chernclasses}
The Chern classes $c_1$ and $c_2$ can be computed easily by
restriction to a divisor. In general, if $Y$ is a smooth
divisor on $X$, the exact sequence
$$
0\longrightarrow {\cT}_Y\longrightarrow {\cT_X}_{|Y}\longrightarrow
\cN_{Y/X}\longrightarrow 0
$$
gives the following relations in $A^*(Y)$:
$$
c_1(X)_{|Y}=c_1(Y)+Y^2,\quad
c_2(X)_{|Y}=c_1(Y)\cdot Y^2+c_2(Y)
,\quad c_3(X)_{|Y}=c_2(Y)\cdot Y^2+c_3(Y)
$$
where $Y^2=c_1(N_{Y/X})\in A^1(Y)$, the first relation is
equivalent to the adjunction formula.

\subsection{Proposition.}\label{chernbo}
The Chern classes of a boundary divisor
$B_0$ are:
$$
c_1(B_0)=(2\hB_b+\hC_b)/5,\qquad
c_2(B_0)=(4\hB_b^2-36\hC_b^2)/25,\qquad
c_3(B_0)=34.
$$

\ts
To compute $c_2(B_0)$ we consider the divisor $Q\cong (\PP^1)^2$
which lies over a cusp in the Segre threefold.
We can identify $Q$ with the divisor $D_1=B_0\cap C_0$ as in
Proposition \ref{intcusp}. Only six of the nine boundary divisors,
besides $B_0$,
which meet $C_0$ intersect $Q$, they intersect it in
$3D_1D_2+3D_1D_3=3\Delta$ where $\Delta$ is the diagonal of $(\PP^1)^2$.
Thus $\hB_b$ restricts to $3\Delta$ on $Q$. Since $\hC_b$ is the sum of
the $10$ disjoint cusp divisors on $B_0$, one of which is $Q$, its
restriction to $Q$ is $Q^2=N_{Q/B_0}$. Since $K_Q=-2\Delta$ the
adjunction formula gives
$$
-2\Delta=K_{B_0|Q}+Q^2=\left((-2\hB_b-\hC_b)/5\right)_{|Q}+Q^2=
(-6\Delta+4Q^2)/5
$$
hence $Q^2=-\Delta$. Furthermore, since
$\hB_b^2\hC_b=10\cdot 18$, $\hC_b^3=20$ and $\hC_b$ is the sum
of $10$ disjoint copies of $Q$,
we get $\hB_b^2Q=18$, $\hC_b^2Q=2$.

Now we write $c_2(B_0)=n\hB_b^2+m\hC_b^2$ and
use the relation \ref{chernclasses} for $Q\subset B_0$:
$$
(n\hB_b^2+m\hC_b^2)_{|Q}=c_1(Q)Q^2+c_2(Q)=(2\Delta)(-\Delta)+4=0,
$$
hence $18n+2m=0$. On any smooth threefold we have
$c_1c_2/24=\chi(\cO)$.
Since $\chi(\cO)=1$ for $B_0$ this implies
$$
24=\mbox{$\frac{1}{5}$}(2\hB_b+\hC_b)(n\hB_b^2+m\hC_b^2)=
-30n-20m.
$$
From these equations one finds $n=4/25$ and $m=-36/25$.

Recall that $B_0$ is obtained as a blow up
in two different ways.
One is described in section \ref{bbu}, the other
in the proof of Proposition \ref{intB0}.
In \cite{FultonI}, Thm.\ 15.4 the Chern
classes of the blow up of a smooth threefold $Y$ in a point $p$
or in a smooth curve $X$ are given.
Applying these formulas allows one to compute,
in two different ways, $c_3(B_0)$ and we found the same result,
$c_3(B_0)=34$, in both ways. This method also allowed us to verify
the formula for $c_2(B_0)$ determined above.
\qed

\subsection{Theorem.} The Chern classes of the cross ratio variety
$\cC$ are:
$$
c_1(\cC)=(\hB-\hC)/4,\qquad
c_2(\cC)=(\hB^2-9\hC^2)/8,\qquad
c_3(\cC)=\mbox{$\frac{13}{3\cdot 96}$}\hB^3+
\mbox{$\frac{181}{96}$}\hC^3,
$$
and $c_4(\cC)=271$.

\ts
Since $c_1(\cC)=-K_{\cC}$, the first result follows
from Theorem \ref{chowC}.
The expression for $c_2$ follows from the formula for $c_2(X)_{|Y}$
 in \ref{chernclasses} applied to $B_0\subset \cC$ and the previous
proposition. Similarly, since
the Chern class $c_3$ is invariant under $W(E_6)$ it
can be written as $a\hB^3+b\hC^3$ for certain $a,\,b\in\QQ$.
Restricting to $B_0$ and to a cusp divisor
$C_0$ and using $c_3(B_0)=34$, $c_3(C_0)=8$
(which follows from $C_0\cong(\PP^1)^3$)
we find two equations for $a$ and $b$ from which $c_3$
can be determined. The class $c_4$ can be computed from the other
$c_i$ and the relation $td(\cT)_4=1$.
\qed

\section{Applications}

\subsection{The CAF-linear system.}
In \cite{Freitag}, Freitag determined the Hilbert function of the
graded ring $\CC[Y_0,\ldots,Y_9]/\cI$, where $\cI$ is the homogeneous
ideal defining $\cM\subset\PP^9$. The Hilbert function
coincides with the map $n\mapsto\chi(\cO_\cC(nH))$,
see the theorem below.
This suggests that
that the multiplication maps
$Sym^nH^0(H)\rightarrow H^0(nH)$ are surjective
and that $H^i(nH)=0$ for $i,n>0$.

\subsection{Theorem.}\label{poli}
Let $H=(\hB+3\hC)/4$ be the pull-back of the hyperplane section
of $\cC\rightarrow \PP^9$. Then we have:
$$
\chi(\cO_\cC(nH))=
\mbox{$\frac{9}{8}$}n^4+
\mbox{$\frac{9}{4}$}n^3+
\mbox{$\frac{27}{8}$}n^2+
\mbox{$\frac{9}{4}$}n+1.
$$

\ts This follows by direct computation from the Riemann-Roch formula.
\qed

\subsection{Tritangent divisors.}\label{tritangent}
We haven't used the tritangent divisors to obtain the results above.
However, applying the same methods we can describe the Chow groups
of a tritangent divisor. We will denote the tritangent divisor in $\cC$
obtained from $\PP^3_{\rm w}$ in $\tT'$ by $T_0$. The strict transform
of the hyperplane class in $\PP^3_{\rm w}$ will be denoted by
$H_{\rm w}\in A^1(T_0)$. There is an identification
$\PP^3_{\rm w}\cong \PP(Q(F_4)\otimes_\ZZ\CC)$, where $Q(F_4)$ is
the root lattice of the root system
system of type $F_4$ such that the action of $W(F_4)$ on $Q(F_4)$
induces the biregular action of $W(F_4)\;(\subset W(E_6))$ on $T_0$
(cf.\ \cite{vG}, 6.4, 6.5).
The root system $F_4$ has $24$ short roots,
these are the elements of $S$,
and 24 long roots, the elements of $R$, these form a subroot system
of type $D_4$ in $F_4$.
The 12 points in $\PP(Q(F_4)\otimes_\ZZ\CC)$
defined by the short roots of $F_4$ are blown up to divisors
in $T_0$.
The sum of these $12$ divisors in $T_0$
will be denoted by $\hB_{t,e}$.

\subsection{Proposition.}\label{chowtrit}
The Chow groups of a tritangent divisor $T_0$ are:
$$
A^1(T_0)\cong A^2(T_0)\cong\ZZ^{29}.
$$
Moreover we have:
$$
A^1(T_0)^{W(F_4)}\cong A^2(T_0)^{W(F_4)}\cong\ZZ^3
$$
and the classes
$H_{\rm w}$, $\hB_t:=\hB_{|T_0}$ and $\hC_t:=\hC_{|T_0}$
are a basis of $A^1(T_0)^{W(F_4)}_\QQ$.

The chern classes of $T_0$ are:
$$
-c_1(T_0)=K_{T_0}=2H_{\rm w} -\frac{1}{2}\hB_t-\frac{1}{2}\hC_t,\qquad
c_2(T_0)=6H_{\rm w}^2 -2H_{\rm w}\hC_t-\hC_t^2,
\qquad c_3(T_0)=92.
$$

The hyperplane class of the CAF-map restricts to
$$
3H_{\rm w}-\hB_{t,e}.
$$

\ts
The 12 curves $C_i'\subset \tT'$
which are blown up in $\pi'$ intersect
the projective space $\PP^3_{\rm w}$
in $12$ points, these are the points defined
by $S$. This adds $12$ to
the rank of $A^1$. Next in $\tT''$ the strict transform
of $\PP^3_{\rm w}$ meets the 16 $S_i''$'s
in 16 smooth rational curves (these curves are the strict
transforms of $16$ lines in $\PP^3_{\rm w}$, each line contains three
of the points, on each point there are 4 lines).
Thus for the strict transform $\hat{E}$ of
$\PP^3_{\rm w}$ in $\hT$ we get rk$A^1(\hat{E})=1+12+16=29$.

The intersection of $\hat{E}$ and a $\hat{C}_i\cong \PP^1\times V$
is a copy of $V$ (cf.\ \cite{Naruki}, p.\ 22).
Since the contraction
$r:\hT\rightarrow \cC$ induces the projection
$\PP^1\times V\rightarrow V$, the map
$r$ induces an isomorphism $\hat{E}\cong T_0$.
Thus also rk$A^1(T_0)=29$.
The rank of $A^2$ follows
since $T_0$ is easily seen to be smooth by construction.
The invariant classes are also obvious.

A tritangent divisor has non-empty intersection with
each boundary divisor.
The $12$ divisors in $\tT'$ obtained from the subtori $\alpha=1$ in $T$,
with $\alpha$ a root of $D_4$,
intersect $\PP^3_{\rm w}$ in linear subspaces
defined by the same root of $F_4$ (\cite{vG}, proof of Thm.\ 6.5).
We denote by $\hB_{t,i}$ the sum of the twelve divisors in $T_0$
obtained in this way.
The other $36-12=2\cdot 12$ boundary divisors in $\cC$ are obtained
from the $V(\tau)$ where $\tau\in S$, these are paired by
$\tau\leftrightarrow -\tau$.
The image $C_j$ in $\tT$ of each of the twelve curves $C_j'$ in $\tT'$
has non-empty intersection with exactly two such boundary divisors,
the two are in a pair and each intersection is transversal in one point.
We denote the pair defined by $C_j$ by $\pm \tau_j$.
The divisors $\hat{T}$ and $\widehat{V(\pm\tau_j)}$ are still disjoint,
but each meets the divisor
$\hat{C}_j\cong C_j\times V$ in a copy of $V$.
In the final blow down $r:\hat{T}\rightarrow \cC$,
$\hat{C}_j$ is contracted onto $V$, hence this surface is
the intersection of any two of $T_0$, $B_{\tau_j}$, and $B_{-\tau_{j}}$
where $B_\tau$ is the boundary divisor defined by $V(\tau)$.
We denote the sum of these 12 surfaces by $\hB_{t,e}$, then
$$
\hB_t=\hB_{t,i}+2\hB_{t,e}.
$$

Using the divisor of $\lambda-1$ (cf.\  \ref{divlambda})
restricted to $T_0$ we find that $T_0^2$ is is the sum of $(-1)$ times
one of the twelve divisors in $\hB_{t,i}$
(it is $D^1_\lambda\cap T_0)$), minus the sum of $4$ divisors in $\hC_t$
and plus the sum of $6$ divisors from $\hB_{t,e}$.
Symmetrizing w.r.t.\ $W(F_4)$ (so we replace each divisor in
$\hB_{t,i}$ by $(1/12)\hB_{t,i}$ etc.) we find:
$$
T_0^2=\frac{-1}{12}\hB_{t,i}+\frac{6}{12}\hB_{t,e}-\frac{4}{16}\hC_t=
\frac{-1}{12}\hB_{t,i}+\frac{1}{2}\hB_{t,e}-\frac{1}{4}\hC_t.
$$

The canonical divisor of $T_0$ can be obtained from the
adjunction formula for $T_0\subset\cC$:
$$
K_{T_0}={K_{\cC}}_{|T_0}+T_0^2=\frac{1}{4}(-\hB_{t,i}-2\hB_{t,e}+\hC_t)
-\frac{1}{12}\hB_{t,i}+\frac{1}{2}\hB_{t,e}-\frac{1}{4}\hC_t=
\frac{-1}{3}\hB_{t,i}.
$$
The canonical divisor of $T_0$ can also be computed from
the construction of $T_0$ given above:
$$
K_{T_0}=-4H_{\rm w}+2\hB_{t,e}+\hC_t.
$$
Thus we find:
$$
\hB_{t,i}=12H_{\rm w}-6\hB_{t,e}-3\hC_t.
$$
The $12H_{\rm w}$ in this formula arises from the fact that each of
the twelve divisors in $\hB_{t,i}$ is the strict transform of
plane in $\PP^3_{\rm w}$. Since $\hB_t=\hB_{t,i}+2\hB_{t,e}$,
this relation implies $4\hB_{t,e}=12H_{\rm w}-\hB_t-3\hC_t$
and thus we can eliminate $\hB_{t,e}$ from the
last formula for $K_{T_0}$ to obtain the formula in the theorem.
The other Chern classes can be computed by restricting the Chern
classes of $\cC$ or directly by
using the blow up formulas as in the
proof of \ref{chernbo}.

The restriction of the CAF-hyperplane class to $T_0$ is:
$$
\mbox{$\frac{1}{4}$}(\hB+3\hC)_{|T_0}=
\frac{1}{4}(\hB_{t,i}+2\hB_{t,e}+3\hC_t)=3H_{\rm w}-\hB_{t,e}
$$
which agrees with the fact that the CAF-map, when restricted to
$\PP^3_{\rm w}$, is given by the cubics which vanish on the
12 points which are blown up to give $\hB_{t,e}$ (\cite{vG}, Thm 6.5).
\qed

\subsection{Remark.}
The tritangent divisor $T_0$ also has non-empty intersection with
the other $45-1=12+2\cdot 16$ tritangent divisors.
It is verified in \cite{vG}, 6.7 that $2\cdot 16$ tritangent divisors
intersect $T_0$ in the strict transform of 16 quadrics in
$\PP^3_{\rm w}$.

The other 12 tritangent divisors
intersect $\PP^3_{\rm w}$ in the $12$ linear subspaces defined by the
short roots, that is by a pair $\pm \tau\in S$.
The involutions corresponding to these 12 divisors are just
the reflections in the corresponding short roots in
$W(F_4)\subset W(E_6)$ (\cite{vG}, proof of Thm 6.5).
The involution defined by $\tau$
permutes the boundary divisors $B_{\pm \tau}$
(defined by $V(\pm\tau)\subset\tT$) in $\cC$. Thus we conclude that
the surface $V=T_0\cap B_{\tau}\cap B_{-\tau}$ must coincide with the
divisor in $T_0$ defined by the plane in $\PP^3_{\rm w}$ defined by
the short root $\tau$. In particular, $V$ lies in two tritangent
as well as in two boundary divisors.

The fact that there are surfaces in the intersection
of pairs of boundary and
tritangent divisors can also easily be seen using the ball quotient
description of $\cC$. Following \cite{AF},
the ball ${\cal B}$ is viewed as a subvariety of $\CC^5$:
$$
{\cal B}=\{(1,z_1,z_2,z_3,z_4)\in\CC^5:\; 1-|z_1|^2-\ldots-|z_4|^2>0\,\}
$$
and we consider the hermitian form $(z,w)=z_0\bar{w}_0-z_1\bar{w}_1-
\ldots -z_4\bar{w}_4$ on $\CC^5$.
Let $e_1=(0,1,0,0,0)$, $e_2=(0,0,1,0,0)$, then the subvarieties
of the ball defined by $(z,e_i)=0$ map onto tritangent divisors $T_i$
under the map
$$
{\cal B}\longrightarrow {\cal B}/\Gamma\cong\cC-\{\text{cusp divisors}\}.
$$
The subvarieties defined by $w_1=e_1+e_2$ and $w_2=e_1-e_2$
map onto boundary divisors $T_j$. Obviously, the intersection of these
four divisors is a surface.

\section{Appendix}\label{app}

\subsection{The divisor of $\lambda$.}\label{divlambda}
The divisor of the rational function $\lambda$, as well as $\lambda-1$,
are used to obtain relations in the various Chow groups considered in
this paper.  In the
toric variety $\tT$ we have
$$
(\lambda)=\sum_{\tau\in S\cup R}
n_\tau V(\tau) \qquad{\rm with}\quad
n_\tau=(\epsilon_1-\epsilon_2,\tau)
$$
(standard inner product on $\RR^4$).
Obviously, we have $n_{-\tau}=-n_{\tau}$, so it suffices to list
the $\tau\in R\cap S$ with $n_{\tau}>0$ which is easily done.
$$
n_\tau>0\quad{\rm iff}\quad S_{\lambda,0}\cup R_{\lambda,0}\cup
\{\epsilon_1-\epsilon_2\},
$$
where
$$
S_{\lambda,0}=\{\epsilon_1,-\epsilon_2,
(\epsilon_1-\epsilon_2+\epsilon_3+\epsilon_4)/2,
(\epsilon_1-\epsilon_2-\epsilon_3-\epsilon_4)/2,
(\epsilon_1-\epsilon_2+\epsilon_3-\epsilon_4)/2,
(\epsilon_1-\epsilon_2-\epsilon_3+\epsilon_4)/2
\},
$$
here $\tau\in S$ so these give boundary
divisors,
$$
R_{\lambda,0}=\{
\epsilon_1+\epsilon_3,\epsilon_1-\epsilon_3,
\epsilon_1+\epsilon_4,\epsilon_1-\epsilon_4,
-\epsilon_2+\epsilon_3,-\epsilon_2-\epsilon_3,
-\epsilon_2+\epsilon_4,-\epsilon_2-\epsilon_4
\}
$$
here $\tau\in R$ so these give cusp
divisors, and finally $n_\tau=2$ if $\tau=\epsilon_1-\epsilon_2$,
the divisor $V(\epsilon_1-\epsilon_2)$ is also a cusp divisor.

To find the divisor of the rational function $\lambda$ on
$\cC$, we observe that a divisor $V(\tau)$ does not
intersect $T$, hence it does not contain $e\in T$ (which gets blown
up in $\tT'$) nor any of the 12 curves in $\tT'$
which are blown up in $\tT''$ nor any of the 16 surfaces in $\tT''$
which are blown up in $\tT$. Thus the divisor of $\lambda$ in $\hT$ is
simply the sum of the strict transforms of the divisors $V(\tau)\subset
\tT$ with the same multiplicities.
Finally, the direct image under $r$
of this divisor in $\cC$ is the divisor
of $\lambda$ on $\cC$.

Using the tables from \ref{bd} below one finds that
$$
S_{\lambda,0}\!=\!\{h_{345},h_{13},h_{136},h_{26},h_{246},h_{256}\}
\subset\{h_{345},h_{13},h_{145},h_{136},h_{26},h_{123},
h_{246},h_{256},h_{45}\}\!=\![13.45.26],
$$
that is, this set of $6$ positive roots can be completed to a set
three orthogonal $A_2$'s (see \ref{cd}). Similarly one finds
$$
S_{\lambda,\infty}\!=\!\{h_{12},h_{245},h_{36},h_{126},h_{356},h_{346}\}
\!\subset\!\{h_{12},h_{245},h_{145},h_{36},h_{126},h_{123},
h_{356},h_{346},h_{45} \}\!=\![12.45.36],
$$
so, using the {\em same} three roots $h_{145},h_{123},h_{45}$,
also this set can be completed
to a set of three orthogonal $A_2$'s.
These sets of nine roots
correspond to the cusp divisors $D_v$ where $v=(1,-1,1,0,0)\in N_0$
(with
$D_v=V(\epsilon_1-\epsilon_2)$) and $w=(1,1,-1,0,0)$ (with
$D_w=V(-\epsilon_1+\epsilon_2)$) respectively. For $z\in N_0$
we define $\{z\}$ to be the set of nine (positive) roots in the
corresponding $A_2^3$ and we put
$$
B_z=\sum_{t\in \{z\}} D_t.
$$

Next we consider the cusp divisors.
Using the table \ref{cd}
we find that the eight $\tau\in R_{\lambda,0}$ correspond to the
lines generated by the eight vectors
$(1,-1,0,\pm1,0)$, $(1,-1,0,0,\pm1)$,
$(1,0,1,\pm1,0)$, $(1,0,1,0,\pm1)$. These vectors are all
perpendicular to $w$, similarly the $\tau\in R_{\lambda,\infty}$
correspond to the vectors $(1,1,0,\pm1,0)$, $(1,1,0,0,\pm1)$,
$(1,0,-1,\pm1,0)$, $(1,0,-1,0,\pm1)$
which are all in $v^\perp$. The linear space $v^\perp\cong \FF_3^4$
contains $13$ isotropic lines, $8$ of these are spanned by the
vectors listed, another $4$ are spanned by $(0,1,1,\pm1,0)$
and $(0,1,1,0,\pm1)$
(note that these vectors are also in $w^\perp$) and finally there
is the line spanned by $v$ itself. We write
$$
C_{v^\perp}=\sum_{z\in v^\perp} D_z
$$
for the sum of these 13 cusp divisors.
Putting all this together, we have
$(\lambda)=B_v+C_{w^\perp}+3D_v-(B_w+C_{v^\perp}+3D_w)$.

As a byproduct of our labors we get:
$$
E_v=E_w,\quad(\in A^1(\cC))\qquad {\rm with}\quad
E_v=B_v-C_{v^\perp}+3D_v.
$$
The group $W(E_6)$ permutes the $v\in N_0$
and it easily follows that any two $E_v$'s are linearly equivalent,
hence $E_v$ is an invariant class. Symmetrizing $E_v$ gives:
$$
E_v=\frac{9}{36}\hB+\frac{(-13+3)}{40}\hC=(\hB-\hC)/4=-K_\cC.
$$

\subsection{The divisor of $\lambda-1$.}\label{divlambda1}
The rational function $\lambda-1$ on has the same poles
as $\lambda$.
Its divisor of zeroes $(\lambda-1)$ in $\tT$
is the closure in $\tT$
of the subtorus defined by $\lambda=1$ in $T$,
we denote it by $D_\lambda^1$.
Since $\lambda(e)=1$, the divisor $(\lambda-1)_0$ in $\tT'$,
the blow up $\tT$ of $e$, is given by the sum of the
strict transform of $D^1_\lambda$ and the exceptional fibre
$\PP^3_{\rm w}$. The strict transforms in $\hT$ of the
exceptional divisors for the map
$\pi':\tT''\rightarrow \tT'$ are blown down by $r:\tT\rightarrow \cC$,
hence in the end they do not contribute to $(\lambda-1)_0$ in
$A^1(\cC)$.
The map $\pi''$ blows up 16 surfaces in $\tT'$, these
are the strict transforms of surfaces in $\tT$. By inspection one finds
that 4 of these are contained in $D^1_\lambda$, to be explicit,
the defining equations of these curves are:
$$
\lambda=\mu\nu\rho=1,\qquad
\lambda=\mu\rho=1,\qquad
\lambda=\nu\rho=1,\qquad
\lambda=\rho=1.
$$
Therefore the divisor $(\lambda-1)_0$ in $\hT$ has 6 components
(with multiplicity one), the strict transforms of
$D^1_\lambda\;(\subset\tT)$, of
$\PP^3_{\rm w})\;(\subset \tT')$
and the exceptional divisors over the 4 surfaces.
The divisor $(\lambda-1)_0$ in $\cC$ is the direct image of this divisor,
in particular it is the sum of one boundary divisor (coming from
$D^1_\lambda$), one tritangent divisor (coming from $\PP^3_{\rm w}$)
and 4 cusp divisors (coming from the exceptional divisors over
the surfaces).

\subsection{Tables.} We collect some of the notations used by
various authors for referring to tritangents, boundary and
cusp divisors (cf.\ \cite{Freitag2} and \cite{Sek}).
We also give a $W(E_6)$-equivariant map between the sets of these
divisors and the sets $N_i$ (see \ref{ni}).

\subsection{The boundary divisors.}\label{bd}
The $36$ boundary divisors are parametrized by the $36$ positive roots
of $E_6$. These roots generate the root lattice $Q(E_6)\cong \ZZ^6$.
The quadratic form given by the Cartan matrix
on $Q(E_6)$ has determinant $3$, hence the $\FF_3$-valued
bilinear form it induces on $FF^6_3$ is degenerate. In fact,
$v=\alpha_1-\alpha_3+\alpha_5-\alpha_6=h_{12}-h_{23}+h_{45}-h_{56}\in
Q(R)$ has the property that $(v,\alpha)\in 3\ZZ$
(this is easy to verify on a basis of simple roots of $E_6$
using the Dynkin diagram),
but $v\not\in 3Q(R)$.
More intrinsically,
we note that $v/3$ lies in the weight lattice $P(E_6)$ and that
$P(E_6)=Q(E_6)+\ZZ(v/3)$, hence
$Q(E_6)/(3Q(E_6)+\ZZ v)\cong Q(E_6)/3P(E_6)$ (cf.\ \cite{Bo},
Exercises, VI \S 4, n.\ 2).
On $\FF_3^5\cong Q(E_6)/(3P(E_6)$ we get a non-degenerate
$\FF_3$-valued quadratic form $q$ and we have an induced homomorphism
$W(E_6)\rightarrow O(\FF_3^5,q)$.

Under the quotient map
$$
\pi:Q(R)\longrightarrow Q(E_6)/(3Q(E_6)+\ZZ v)\cong \FF_3^5
$$
the roots of $E_6$ map to elements of $N_{-1}$.
We fix $\pi$ by defining the images of the basis of simple roots
$h_{12}\,,h_{123},\,h_{23},\,h_{34},h_{45},\, h_{56}$:
$$
\pi(h_{i,i+1})=f_i-f_{i+1}\qquad(1\leq i\leq 4),\qquad
\pi(h_{123})=(0,0,0,1,1),
$$
here $f_i$ is the $i$-th standard basis vector of $\FF_3^5$,
and in view of $\pi(v)=0$ we get
$$
\pi(h_{56})=h_{12}-h_{23}+h_{45}=(1,1,1,1,-1).
$$
The first table lists the roots of $D_4\subset E_6$, the corresponding
characters of $T$, the name of the root in $E_6$ and its image in
$\PP(\FF_3^5)$. The last three tables list the elements in $S$, which
define boundary divisors in $\cC$, the corresponding
positive root in $E_6$ and its image in $\PP(\FF_3^5)$.
Note that $(\pm\pm\pm\pm)$ stands for the vector
$(\pm\epsilon_1\pm\epsilon_2\pm\epsilon_3\pm\epsilon_4)/2$.

\[
\begin{tabular}{cc}

\begin{tabular}[c]{|rrrr|}
\hline
$\mathrm{e}_{3}+\mathrm{e}_{4} $&$\quad\mu $&$\quad h_{123}$&$ \quad
(0,0,0,1,1)$\\
\hline
$\mathrm{e}_{3}-\mathrm{e}_{4} $&$ \nu $&$h_{45}$&$ (0,0,0,1,-1)$\\
\hline
$\mathrm{e}_{2}+\mathrm{e}_{4} $&$\mu\rho $&$h_{124}$&$(0,0,1,0,1)
$\\
\hline
$\mathrm{e}_{2}-\mathrm{e}_{4}$&$ \nu\rho $&$h_{35} $&$ (0,0,1,0,-1)$\\
\hline
$\mathrm{e}_{2}+\mathrm{e}_{3}$&$\mu\nu\rho $&$h_{125} $&$ (0,0,1,1,0)$\\
\hline
$\mathrm{e}_{2}-\mathrm{e}_{3} $&$ \rho $&$  h_{34} $&$ (0,0,1,-1,0)$\\
\hline
$\mathrm{e}_{1}+\mathrm{e}_{4} $&$ \lambda\mu\rho $&
      $  h_{134} $&$ (0,1,0,0,1)$\\
\hline
$\mathrm{e}_{1}-\mathrm{e}_{4} $&$ \lambda\nu\rho $&
     $  h_{25} $&$ (0,1,0,0,-1)$\\
\hline
$\mathrm{e}_{1}+\mathrm{e}_{3} $&$ \lambda\mu\nu\rho$&
   $  h_{135} $&$ (0,1,0,1,0)$\\
\hline
$\mathrm{e}_{1}-\mathrm{e}_{3} $&$ \lambda\rho $&
   $  h_{24} $&$ (0,1,0,-1,0)$\\
\hline
$\mathrm{e}_{1}+\mathrm{e}_{2}$&$\lambda\mu\nu\rho^{2}$&
 $h_{145}$&$(0,1,1,0,0)$\\
\hline
$\mathrm{e}_{1}-\mathrm{e}_{2}$&$\lambda$&$h_{23}$&$(0,1,-1,0,0)$\\
\hline
\end{tabular}

&

\begin{tabular}
[c]{|rrr|}
\hline
$\epsilon_{4} $&$ \quad h_{234} $&$\quad (1,0,0,0,1)$\\
\hline
$-\epsilon_{4} $&$ h_{15} $&$ (1,0,0,0,-1)$\\
\hline
$\epsilon_{3} $&$ h_{235} $&$ (1,0,0,1,0)$\\
\hline
$-\epsilon_{3} $&$ h_{14} $&$ (1,0,0,-1,0)$\\
\hline
$\epsilon_{2} $&$ h_{245} $&$ (1,0,1,0,0)$\\
\hline
$-\epsilon_{2} $&$ h_{13} $&$ (1,0,-1,0,0)$\\
\hline
$\epsilon_{1} $&$ h_{345} $&$ (1,1,0,0,0)$\\
\hline
$-\epsilon_{1} $&$ h_{12} $&$ (1,-1,0,0,0)$\\
\hline
\end{tabular}

\end{tabular}
\]

\[
\begin{tabular}{cc}

\begin{tabular}[t]{|rrr|}
\hline
$(----) $&$ h $&$ (1,1,1,1,1)$\\
\hline
$(---+) $&$ h_{56} $&$ (1,1,1,1,-1)$\\
\hline
$(--+-) $&$ h_{46} $&$ (1,1,1,-1,1)$\\
\hline
$(--++) $&$ h_{456} $&$ (1,1,1,-1,-1)$\\
\hline
$(-+--) $&$ h_{36} $&$ (1,1,-1,1,1)$\\
\hline
$(-+-+) $&$ h_{356} $&$ (1,1,-1,1,-1)$\\
\hline
$(-++-) $&$ h_{346} $&$ (1,1,-1,-1,1)$\\
\hline
$(-+++) $&$ h_{126} $&$ (1,1,-1,-1,-1)$\\
\hline
\end{tabular}

&

\begin{tabular}[t]{|rrr|}
\hline
$(+---) $&$ h_{26} $&$(1,-1,1,1,1)$\\
\hline
$(+--+) $&$ h_{256} $&$ (1,-1,1,1,-1)$\\
\hline
$(+-+-) $&$ h_{246} $&$ (1,-1,1,-1,1)$\\
\hline
$(+-++) $&$ h_{136} $&$ (1,-1,1,-1,-1)$\\
\hline
$(++--) $&$ h_{236} $&$ (1,-1,-1,1,1)$\\
\hline
$(++-+) $&$ h_{146} $&$ (1,-1,-1,1,-1)$\\
\hline
$(+++-) $&$ h_{156} $&$ (1,-1,-1,-1,1)$\\
\hline
$(++++) $&$ h_{16} $&$ \quad (1,-1,-1,-1,-1)$\\
\hline
\end{tabular}

\end{tabular}
\]

\subsection{Tritangents.}\label{td}
The tritangents were labelled by Cayley, later a more comprehensible
notation was introduced by Schl\"afli. The $27$ lines on a smooth cubic
surface are denoted by $a_1,\ldots,a_6$ (these are $6$ skew lines which
can thus be blown down to points $p_i$),
$b_1,\ldots,b_6$ ($b_i$ maps to the
conic on the $6$ points except for $p_i$) and $c_{ij}$,
$1\leq i<j\leq 6$ ($c_{ij}$ maps to the line $\langle p_i,p_j\rangle$).
The 45 tritangent planes are denoted by
$$
(ij)=\langle a_i,b_j,c_{ij}\rangle,\qquad
(ij.kl.mn)=\langle c_{ij},c_{kl},c_{mn}\rangle,
$$
here $i\neq j$ and $(ij)\neq (ji)$,
also $\{i,\ldots,n\}=\{1,\ldots,6\}$ and $i<j$, $k<l$, $m<n$.
A correspondence between Cayley's labels and those of Schl\"afli was
given in \cite{Sek} and is copied here.

Each tritangent plane is determined by a $D_4\subset E_6$
(cf.\ \cite{vG}, 1.8). The span of the image of such a
$D_4$ in $\FF_3^5$ has codimension one (the image of the `standard'
$D_4$ spanned by $h_{123},\,h_{23},\,h_{34},\,h_{45}$ has this property,
the other $D_4$'s are in $W(E_6)$-orbit of this one), hence it is the
perpendicular of a one dimensional subspace. In this way
each tritangent defines a point in
$\PP(\FF_3^5)$, one checks that this point lies in $N_1$ for the
standard $D_4$, and thus it lies in $N_1$ for all tritangents.
Alternatively, one can use \cite{AF} to find this correspondence.

Another way to describe this correspondence was found by
Matsumoto and Terasoma \cite{MT}, Prop.\ 3.4. They identify
$\FF_3^5$ with a subgroup of the $3$-torsion points of the intermediate
Jacobian of a cubic threefold associated to a cubic surface.

\[
\begin{tabular}[c]{cc}

\begin{tabular}[c]{|rrr|}
\hline
$\xi $&$  {(56)} $&$ (0,0,0,0,1)$\\
\hline
$x $&$ {(46)} $&$ (0,0,0,1,0)$\\
\hline
$\mathrm{x} $&$ {(36)} $&$ (0,0,1,0,0)$\\
\hline
$\overline{\mathrm{x}} $&$ {(26)} $&$(0,1,0,0,0)$\\
\hline
$\mathrm{w} $&$ {(16)} $&$ (1,0,0,0,0)$\\
\hline
$\overline{\mathrm{y}} $&$ {(61)}$&$ (0,1,1,1,1)$\\
\hline
$\overline{\mathrm{z}}$&$ {(15)}$&$ (0,1,1,1,-1)$\\
\hline
$\mathrm{z} $&$ {(14)} $&$ (0,1,1,-1,1)$\\
\hline
$\mathrm{y} $&$ {(16.23.45)} $&$ (0,1,1,-1,-1)$\\
\hline
$z $&$ {(13)} $&$ (0,1,-1,1,1)$\\
\hline
$y $&$ {(16.24.35)} $&$ (0,1,-1,1,-1)$\\
\hline
$\eta $&$  {(16.25.34)} $&$ (0,1,-1,-1,1)$\\
\hline
$\zeta $&$  {(12)} $&$ (0,1,-1,-1,-1)$\\
\hline
$\overline\mathrm{{r}}_{1} $&$ {(62)} $&$ (1,0,1,1,1)$\\
\hline
$\overline{\mathrm{q}}_{1} $&$ {(25)} $&$ (1,0,1,1,-1)$\\
\hline
$\overline{\mathrm{m}}_{1} $&$ {(24)} $&$ (1,0,1,-1,1)$\\
\hline
$\mathrm{n}_{1} $&$ {(13.26.45)} $&$ (1,0,1,-1,-1)$\\
\hline
$\overline{\mathrm{n}} $&$ {(23)} $&$  (1,0,-1,1,1)$\\
\hline
$m $&$\quad {(14.26.35)} $&$ (1,0,-1,1,-1)$\\
\hline
$\overline{\mathrm{q}} $&$ {(15.26.34)} $&$ (1,0,-1,-1,1)$\\
\hline
$\overline{\mathrm{r}} $&$ {(21)} $&$ (1,0,-1,-1,-1)$\\
\hline
\end{tabular}

&

\begin{tabular}[c]{|rrr|}\hline
$\overline{\mathrm{n}}_{1} $&$ {(63)} $&$ (1,1,0,1,1)$\\
\hline
$\mathrm{m}_{1} $&$ {(35)}$ & $(1,1,0,1,-1)$\\
\hline
$\mathrm{q}_{1} $&$ {(34)}$ & $(1,1,0,-1,1)$\\
\hline
$\mathrm{r}_{1}$ &$ {(12.36.45)}$ & $(1,1,0,-1,-1)$\\
\hline
$\mathrm{r} $&$ {(32)}$ & $(1,-1,0,1,1)$\\
\hline
$\mathrm{q} $&$ {(14.25.36)}$ & $(1,-1,0,1,-1)$\\
\hline
$\overline{\mathrm{m}} $&$ {(15.24.36)}$ & $(1,-1,0,-1,1)$\\
\hline
$\mathrm{n} $&$ {(31)}$ & $(1,-1,0,-1,-1)$\\
\hline
$\overline{\mathrm{l}}^{\phantom{1}} $&$ {(64)}$ &$(1,1,1,0,1)$\\
\hline
$\mathrm{l} $&$ {(45)}$ &$(1,1,1,0,-1)$\\
\hline
$\overline{\mathrm{g}} $&$ {(43)}$ &$(1,1,-1,0,1)$\\
\hline
$\overline{\mathrm{h}} $&$ {(12.35.46)}$ &$(1,1,-1,0,-1)$\\
\hline
$\mathrm{h} $&$ {(42)}$ &$(1,-1,1,0,1)$\\
\hline
$\mathrm{g} $&$ {(13.25.46)}$ & $(1,-1,1,0,-1)$\\
\hline
$\overline{\mathrm{l}}_{1} $&$ {(15.23.46)}$ & $(1,-1,-1,0,1)$\\
\hline
$\mathrm{l}_{1} $&$ {(41)}$ & $(1,-1,-1,0,-1)$\\
\hline
$\overline{\mathrm{p}} $&$ {(65)}$ & $(1,1,1,1,0)$\\
\hline
$\mathrm{p} $&$ {(54)}$ & $(1,1,1,-1,0)$\\
\hline
$\mathrm{f} $&$ {(53)}$ & $(1,1,-1,1,0)$\\
\hline
$\theta $&$ {(12.34.56)}$ & $(1,1,-1,-1,0)$\\
\hline
$\overline{\theta} $&$ {(52)}$ & $(1,-1,1,1,0)$\\
\hline
$\overline{\mathrm{f}}^{\phantom{1}} $&$ {(13.24.56)}$ & $(1,-1,1,-1,0)$\\
\hline
$\mathrm{p}_{1} $&$\quad {(14.23.56)}$ & $(1,-1,-1,1,0)$\\
\hline
$\overline{\mathrm{p}}_{1} $&$ {(51)}$ & $(1,-1,-1,-1,0)$\\
\hline
\end{tabular}

\end{tabular}
\]

\subsection{The cusp divisors.}\label{cd}
The $40$ cusp divisors can be labelled by `triads of trihedral pairs'
(cf.\ \cite{Hunt}, 6.1.1), they correspond also to copies
of three orthogonal $A_2$'s in $E_6$ (cf.\ \cite{Hunt}, 6.1.5.3)
The notation we use is:
$$
[ijk.lmn]=\left[
\begin{array}{ccc}
h_{ij}&h_{jk}&h_{ik}\\
h_{lm}&h_{mn}&h_{ln}\\
h&h_{ijk}&h_{lmn}
\end{array}
\right],
\qquad
[ij.kl.mn]=\left[
\begin{array}{ccc}
h_{ij}&h_{ikl}&h_{jkl}\\
h_{kl}&h_{kmn}&h_{lmn}\\
h_{mn}&h_{nij}&h_{mij}
\end{array}
\right],
$$
here the rows of the matrices are thee positive vectors of $E_6$
which span an $A_2$ and the roots in different rows are perpendicular.

The image under $\pi$ of the span of an
$A_2\subset Q(E_6)$ is a two dimensional
subspace of $\FF_3^5$, its projectivizations thus has $4$ points,
three of which are lines spanned by the positive roots, the other
line is spanned by an element from $N_0$.
The span of the images of three
perpendicular $A_2$'s has codimension one in $\FF_3^5$
(in fact, the span of each of the three contains the same element
$v$ from $N_0$).
Hence this
subspace is the perpendicular of  an element, which in fact is
this $v\in N_0$
(to check this
it suffices to verify this statement for one of the $40$ triples,
and to apply $W(E_6)$ to obtain the result in general).
Conversely, given a vector $v\in N_0$, it is perpendicular
to exactly $9$ elements in $N_{-1}$.

The first table lists the equations of a subtorus of $T$
whose closure in the toric variety $\tT$ is one of the 16 surfaces
$S_i$, the pair of roots in $D_4\subset E_6$
which define the surface are also given. This pair of roots spans an
$A_2\subset E_6$ which is contained in a unique triple of orthogonal
$A_2$'s.
The corresponding cusp divisor in $\cC$ is then labelled as above
as well as by a vector in $N_0$, this vector is normalized by
its first non-zero component being $+1$.
The last table lists a $\tau\in R$, which defines a cusp divisor
in $\cC$, a label for the corresponding triple of orthogonal $A_2$'s
as well as the corresponding normalized vector in $N_0$.

\[
\begin{tabular}
[c]{|r|r|r|r|}\hline
$ \nu = \rho =1 $ & $h_{34 },\,h_{45}$ &
$[126.345]$ & $(0,0,1,1,1)$\\\hline
$ \mu = \rho =1 $ & $h_{34 },\,h_{123}$ &
$[34.12.56] $ & $(0,0,1,1,-1)$\\\hline
$ \mu = \nu\rho =1 $ & $h_{35 },\,h_{123}$ &
$[35.12.46]$ & $(0,0,1,-1,1)$\\\hline
$ \nu = \mu\rho =1 $ & $h_{45 },\,h_{124}$ &
$[45.12.36]$ & $(0,0,-1,1,1)$\\\hline
$ \nu = \lambda\rho =1 $ & $h_{24 },\,h_{45}$ &
$[136.245]$ & $(0,1,0,1,1)$\\\hline
$ \mu = \lambda\rho =1 $ & $h_{24 },\,h_{123}$ &
$[24.13.56]$ & $(0,1,0,1,-1)$\\\hline
$ \mu = \lambda\nu\rho =1$&$h_{123 },\,h_{25}$ &
$[25.13.46]$ & $(0,1,0,-1,1)$\\\hline
$ \nu = \lambda\mu\rho =1 $&$h_{45 },\,h_{134}$ &
$[45.13.26]$ & $(0,1,0,-1,-1)$\\\hline
$ \lambda = \nu\rho =1 $ & $h_{23 },\,h_{35}$ &
$[146.235]$ & $(0,1,1,0,1)$\\\hline
$ \lambda = \mu\rho =1 $ & $h_{23 },\,h_{124}$ &
$[23.14.56]$ & $(0,1,1,0,-1)$\\\hline
$ \lambda\nu\rho = \mu\rho =1 $ & $h_{25},\,
h_{124}$ & $[25.14.56]$ & $(0,1,-1,0,1)$\\\hline
$ \lambda\mu\rho = \nu\rho =1 $ & $h_{35},\,
h_{134}$ & $[35.14.56]$ & $(0,1,-1,0,-1)$\\\hline
$ \lambda = \rho =1 $ & $h_{23 },\,h_{34}$ &
$[234.156]$  & $(0,1,1,1,0)$\\\hline
$ \lambda = \mu\nu\rho =1 $ & $h_{23 },\,h_{125}$ &
$[23.15.46]$ & $(0,1,1,-1,0)$\\\hline
$ \mu\nu\rho = \lambda\rho =1 $ & $h_{24},\,
h_{125}$ & $[24.15.36]$ & $(0,1,-1,1,0)$\\\hline
$ \rho = \lambda\mu\nu\rho =1 $ & $h_{34},\,
h_{135}$ & $[34.15.26]$ & $(0,1,-1,-1,0)$\\\hline
\end{tabular}
\]

\[
\begin{tabular}{cc}
\begin{tabular}
[c]{|rrr|}\hline
$-\epsilon_{3}-\epsilon_{4}$ & $[145.236]$ & $(1,0,0,1,1)$\\
\hline
$-\epsilon_{3}+\epsilon_{4} $ & $\;[14.23.56]$ &$(1,0,0,1,-1)$\\
\hline
$\epsilon_{3}-\epsilon_{4} $ & $[15.23.46]$ & $(1,0,0,-1,1)$\\
\hline
$\epsilon_{3}+\epsilon_{4} $ & $[45.23.16]$ &$(1,0,0,-1,-1)$\\
\hline
$-\epsilon_{2}-\epsilon_{4} $ & $[135.246]$ &$(1,0,1,0,1)$\\
\hline
$-\epsilon_{2}+\epsilon_{4} $ & $[13.24.56]$ &$(1,0,1,0,-1)$\\
\hline
$\epsilon_{2}-\epsilon_{4} $ & $[15.24.36]$ &$(1,0,-1,0,1)$\\
\hline
$\epsilon_{2}+\epsilon_{4} $ & $[35.24.16]$ &$(1,0,-1,0,-1)$\\
\hline
$-\epsilon_{2}-\epsilon_{3} $ & $[134.256]$ &$(1,0,1,1,0)$\\
\hline
$-\epsilon_{2}+\epsilon_{3} $ & $[13.25.46]$ &$(1,0,1,-1,0)$\\
\hline
$\epsilon_{2}-\epsilon_{3} $ & $[14.25.36]$ &$(1,0,-1,1,0)$\\
\hline
$\epsilon_{2}+\epsilon_{3} $ & $[34.25.16]$ &$(1,0,-1,-1,0)$\\
\hline
\end{tabular}

&

\begin{tabular}
[c]{|rrr|}
\hline
$-\epsilon_{1}-\epsilon_{4} $ & $[125.346]$ &$(1,1,0,0,1)$\\
\hline
$-\epsilon_{1}+\epsilon_{4} $ & $\;[12.34.56]$ &$(1,1,0,0,-1)$\\
\hline
$\epsilon_{1}-\epsilon_{4} $ & $[15.34.26]$ & $(1,-1,0,0,1)$\\
\hline
$\epsilon_{1}+\epsilon_{4} $ & $[25.34.16]$ &$(1,-1,0,0,-1)$\\
\hline
$-\epsilon_{1}-\epsilon_{3} $ & $[124.356]$ &$(1,1,0,1,0)$\\
\hline
$-\epsilon_{1}+\epsilon_{3} $ & $[12.35.46]$ & $(1,1,0,-1,0)$\\
\hline
$\epsilon_{1}-\epsilon_{3} $ & $[14.35.26]$ & $(1,-1,0,1,0)$\\
\hline
$\epsilon_{1}+\epsilon_{3} $ & $[24.35.16]$ & $(1,-1,0,-1,0)$\\
\hline
$-\epsilon_{1}-\epsilon_{2} $ & $[123.456]$ & $(1,1,1,0,0)$\\
\hline
$-\epsilon_{1}+\epsilon_{2} $ & $[12.45.36]$ & $(1,1,-1,0,0)$\\
\hline
$\epsilon_{1}-\epsilon_{2} $ & $[13.45.26]$ & $(1,-1,1,0,0)$\\
\hline
$\epsilon_{1}+\epsilon_{2} $ & $[23.45.16]$ & $(1,-1,-1,0,0)$\\
\hline
\end{tabular}
\end{tabular}
\]

\

%\today

\

\end{document}